\documentclass[12pt]{amsart}
\usepackage{latexsym,amsfonts,amsmath,amssymb,amsthm, mathtools,amscd,epsf}
\usepackage[utf8]{inputenc}
\usepackage[dvipsnames]{xcolor}
\usepackage[francais,english]{babel}
\usepackage{tikz}
\usepackage{nicefrac}
\usepackage{stmaryrd}
\usepackage{mathrsfs}
\usepackage[all]{xy}
\usetikzlibrary{matrix,arrows}
\usepackage{lipsum}
\usepackage{graphicx,accents}

\newtheorem{mytheorem}{Theorem}
\newtheorem{theorem}{Theorem}[section]
\newtheorem{lemma}[theorem]{Lemma}
\newtheorem{proposition}[theorem]{Proposition}
\newtheorem{corollary}[theorem]{Corollary}
\newtheorem{definition}[theorem]{Definition}
\newtheorem{remark}[theorem]{Remark}
\newtheorem{example}[theorem]{Example}

\newcommand{\ttor}{( \mathbb{C}^{\star} )^{2}}
\newcommand{\C}{ \mathcal{C}}
\newcommand{\cH}{ \mathcal{H}}
\newcommand{\cL}{ \mathcal{L}}
\newcommand{\cP}{ \mathcal{P}}
\newcommand{\bC}{ \mathbb{C}}
\newcommand{\bZ}{ \mathbb{Z}}

\newcommand{\sw}{ \textsc{w}}
\newcommand{\fn}{\Upsilon}
\newcommand{\scrR}{ \mathscr{R}}

\newcommand{\thb}{\Lambda^{\text{twist}}}

\newcommand{\hc}{K^{\scalebox{0.5}{1/2}}}
\newcommand{\htr}{T^{\scalebox{0.5}{1/2}}}
\newcommand{\sone}{(S^1)^2}
\newcommand{\mgn}{M_{g,n}}
\newcommand{\bmgn}{\overline{M}_{g,n}}

\newcommand{\df}[1]{{\color{NavyBlue} #1}}

\newcommand{\nocontentsline}[3]{}

\newcommand\blfootnote[1]{%
  \begingroup
  \renewcommand\thefootnote{}\footnote{#1}%
  \addtocounter{footnote}{-1}%
  \endgroup
}

\newcommand{\A}{ \mathcal{A}}
\newcommand{\R}{ \mathbb{R}}

\let\Re\relax
\let\Im\relax
\DeclareMathOperator{\SL}{SL}
\DeclareMathOperator{\arcsinh}{arcsinh}
\DeclareMathOperator{\Re}{Re}
\DeclareMathOperator{\Res}{Res}
\DeclareMathOperator{\Im}{Im}
\DeclareMathOperator{\im}{im}
\DeclareMathOperator{\Arg}{Arg}
\DeclareMathOperator{\arccosh}{arccosh}

\DeclareMathOperator{\id}{id}
\DeclareMathOperator{\itr}{int}

\makeindex

\begin{document}

             \title[Harmonic tropical morphisms and approximation]
             {Harmonic tropical morphisms and approximation}

\author[Lionel Lang]{Lionel Lang}
\address{Department of Mathematics, Stockholm University, SE-106 91
Stockholm,         Sweden}
\email {lang@math.su.se }

\maketitle

\begin{abstract}
\textit{Harmonic amoebas} are generalisations of amoebas of algebraic curves immersed in complex tori. Introduced in \cite{Kri}, the consideration of such objects suggests to enlarge the scope of tropical geometry. In the present paper, we introduce the notion of harmonic morphisms from tropical curves to affine spaces and show how these morphisms can be systematically described as limits of families of harmonic amoeba maps on Riemann surfaces. It extends previous results about approximation of tropical curves in affine spaces and provides a different point of view on Mikhalkin's approximation Theorem for regular phase-tropical morphisms, as stated e.g. in \cite{Mikh06}. The results presented here follow from the study of imaginary normalised differentials on families of punctured Riemann surfaces and suggest interesting connections with compactifications of moduli spaces.
\end{abstract} 
 
\blfootnote{The author is supported by the FNS project 140666 and the ERC TROPGEO.\\
The author is grateful to S. Boucksom, E. Brugallé, C. Favre, N. Kalinin, G. Mikhalkin, J. Rau and K. Shaw for many useful discussions and comments. The author is also indebted to an anonymous referee for a careful reading and numerous comments.}

\begin{center}
\begin{Large}
\section{Introduction}
\end{Large}
\end{center}

\subsection{Motivations}

The present work is an attempt to give an appropriate definition of tropical convergence for families of abstract algebraic curves, with a view towards the constructive and enumerative aspects of tropical geometry. In particular, one of the goals of this paper is to provide an alternative proof of Mikhalkin's approximation Theorem \cite[Theorem 1]{Mikh06} on the realisability of phase-tropical curves in $(\bC^\star)^m$.

The link between algebraic and tropical geometry is given by the notion of amoeba introduced in \cite{GKZ}. Recall that the amoeba of an algebraic subvariety $\mathcal{V}\subset (\mathbb{C}^\star)^m$ is the image of $\mathcal{V}$ under the amoeba map 
\[
\begin{array}{rcl}
\A \;\; : \;\; (\mathbb{C}^\star)^m & \rightarrow & \R^m \\
(z_1,\dots,z_m) & \mapsto & \big( \log \vert z_1 \vert, \dots, \log \vert z_m \vert \big)
\end{array}.
\]
Classically, a family of algebraic curves $\left\lbrace \C_t \right\rbrace_{t>1} \subset (\mathbb{C}^\star)^m$ is said to converge to a tropical curve $C \subset \mathbb{R}^m$ if 
\begin{equation}\label{eq:conv}
\lim_{t\rightarrow \infty} \;  \frac{1}{\log (t)} \cdot \A (\C_t) = C 
\end{equation}
with respect to the Hausdorff distance on compact sets. 

In this paper, we aim to understand the relations between the moduli of a family of algebraic curves  converging in the sense of \eqref{eq:conv} and the limiting tropical curve. Conversely, we want to understand how to prescribe the moduli of a family of algebraic curves so that the latter family converges in the sense of \eqref{eq:conv}.

To that purpose, it will be convenient to leave the algebraic setting and consider the more general theory of harmonic amoebas introduced in \cite{Kri}. We will investigate the above problems for a wider class of object that we here call harmonic tropical curves. Eventually, we will come back to the algebraic setting and reprove Mikhalkin's approximation Theorem.

\subsection{Main results}

\subsubsection{Harmonic amoebas and harmonic tropical curves.}$\,$\medskip

For a proper algebraic map $\phi: S \rightarrow (\mathbb{C}^\star)^m$ from a punctured Riemann surface $S$, the composition of the amoeba map $\A$ with $\phi$ is given by integrating the real 1-forms $\Re \big( d \log z_j \big)$ on $S$ where $z_1,\dots,z_m$ are the coordinate functions of $\phi$. The forms $d \log z_j $  are special instances of \df{imaginary normalised differentials} (\df{i.n.d.} for short), namely meromorphic differentials with at worst simple poles at the punctures of $S$ and with purely imaginary periods, see \cite{Kri}. It follows from Riemann's bilinear relations that such a differential is determined by the vector of its residues at the $n$ punctures of $S$. Therefore, a collection of $m$ residue vectors encoded in a real $m \times n$ matrix $R$ defines a collection $\omega_{R,1},\dots, \omega_{R,m}$ of $m$ i.n.d. on $S$ and a map
\[
\begin{array}{rcl}
\mathcal{A}_R \;\;  : \;\;   S & \rightarrow & \mathbb{R}^m \\
s & \mapsto & \Big( \Re\big( \int_p^s \omega_{R,1}\big), \dots, \Re\big( \int_p^s \omega_{R,m}\big)\Big)
\end{array}
\]
up to the choice of a point $p\in S$. Since the coordinates of $\A_R$ are harmonic functions, we call the set $\mathcal{A}_R  (S) \subset \mathbb{R}^m$  the \df{harmonic amoeba}\footnote{The present terminology is due to the author.} of $S$ with respect to $R$. Observe that the map $\A_R$ corresponds to the composition $\A\circ\phi$ for an algebraic map $\phi: S \rightarrow (\mathbb{C}^\star)^m$ if and only if all the periods of the differentials $\omega_{R,1},\dots, \omega_{R,m}$ are integer multiples of $2\pi i$. In such case, the $j^{th}$ coordinate of $\phi$ is given by the exponential of the integral $\int \omega_{R,j}$, see Section \ref{sec:indha} for more details.

In this text, we adapt the above generalisation of amoebas to the case of tropical curves. Following \cite{Mikh05}, tropical curves are graphs with univalent vertices removed and equipped with a complete inner metric (unbounded edges are referred to as leaves). A tropical morphism $\pi:C\rightarrow \R^m$ on a tropical curve $C$ is a continuous map that is affine linear on any edge, with integer slopes and satisfyng the so-called balancing condition, see Section \ref{sec:trop}. Alternatively, tropical morphisms can be  described by integration of exact 1-forms on tropical curves. Following \cite{MZ}, a 1-form $\sw$ on a tropical curve $C$ is equivalent to the data of a current on $C$, seen as an electrical circuit, where a leaf of $C$ is either an electrical source or sink depending whether the residue of $\sw$ at the leaf is positive or negative.

Similarly to i.n.d. on algebraic curves, exact 1-forms on a tropical curve $C$ are in correspondence with residue vectors at the $n$ leaves of $C$, see Proposition \ref{tdiff}. Thus, any $m \times n$ residue matrix $R$ defines a collection $\sw_{R,1},\dots, \sw_{R,m}$ of $m$ exact 1-forms on $C$ and a map
\[
\begin{array}{rcl}
\pi_R \;\;  : \;\;   C & \rightarrow & \mathbb{R}^m \\
q & \mapsto & \Big( \Re\big( \int_p^q \sw_{R,1}\big), \dots, \Re\big( \int_p^q \sw_{R,m}\big)\Big)
\end{array}
\]
up to the choice of a point $p\in C$. We call the set $\pi_R  (C) \subset \mathbb{R}^m$ a \df{harmonic tropical curve}. Observe that the vector of currents induced by $R$ on an edge $e\subset C$ is the slope of $\pi_R(e)$. This vector can alternatively be thought of as a vector of periods of the tropical forms $\sw_{R,1},\dots, \sw_{R,m}$. Therefore, harmonic amoebas (respectively harmonic tropical curves) are honest amoebas (respectively tropical curves) if and only if the periods of the forms involved are in $2\pi i \bZ$ (respectively $\bZ$).

\subsubsection{Tropical convergence and approximation.}$\,$\medskip

To understand the relation between tropical curves and the moduli of families of algebraic curves, we will describe the latter moduli  
using Fenchel-Nielsen coordinates, see Section \ref{sec:hyp} for more details. In this text, we denote $\bmgn$ the (analytic) moduli space of $n$-pointed stable curves of genus $g$ and  $\mgn\subset \bmgn$ the subset of smooth curves.

Recall that any Riemann surface $S\in \mgn$ with $2g+n>2$ admits a unique complete hyperbolic metric, see \cite{SG}. Any collection of $3g-3+n$ hyperbolic geodesics  decomposes $S$ into  $2g-2+n$ pairs of pants. Each pair of pants is uniquely determined by the length of its boundary geodesics (by convention, a puncture is a geodesic of length $0$) 
and two pairs of pants adjacent to a common geodesic are glued to each other in $S$ according to a ``twist" parameter in $S^1$.

Conversely, if $G$ is the graph dual to the above decomposition of $S$ (where pairs of pants and punctures correspond respectively to trivalent and univalent vertices of $G$), we can construct any Riemann surface in  $\mgn$ from a length function $\ell:E(G)\rightarrow\R_{>0}$ and a twist function $\theta:E(G)\rightarrow S^1$ on the set $E(G)$ of edges of $G$ not adjacent to a univalent vertex. We denote the resulting Riemann surface $S(G,\ell,\theta)\in \mgn$ where  the parameters $(\ell,\theta)$ are the aforementioned Fenchel-Nielsen coordinates, relative to $G$.
Observe that the data $(G,\ell)$ is equivalent to the data of a trivalent tropical curve, usually referred to as a simple tropical curve.

\begin{definition}[\textbf{Abstract tropical convergence}]\label{tropconv}$\,$\\
Let $C:=(G,\ell)$ be a simple tropical curve of genus $g$ with $n$ leaves. A family of Riemann surfaces $\left\lbrace S_t \right\rbrace_{t >1} \subset \mgn $ \df{converges to} $C$ if $S_t = S(G,\ell_t, \theta_t)$ for some functions $(\ell_t,\theta_t):E(G) \rightarrow \R_{> 0}\times S^1$ 
and if for any edge $e\in E(C)$, we have
$$ \displaystyle \ell_t(e) \; \underset{t\rightarrow \infty}{\sim} \; \frac{2\pi^2}{\ell(e) \log(t)}.$$
\end{definition}

In the above definition, the family $\left\lbrace S_t \right\rbrace_{t >1}$ converges to the unique stable curve in $\bmgn $ whose dual graph is $G$. The above family defines $3g-3+n$ vanishing cycles in each individual Riemann surface $S_t$. Then, the length function $\ell$ determines the speed of contraction of each vanishing cycle in term of the parameter $t$.

The theorem below establishes the connection between the abstract convergence given in Definition \ref{tropconv} and the notion of convergence given in \eqref{eq:conv}, in the more general context of harmonic amoebas and harmonic tropical curves.

\begin{mytheorem}\label{approxtrop}
Let $C$ be a simple tropical curve of genus $g$ with $n$ leaves and let $R$ be an $m \times n$ matrix of residues. Then, for any sequence $ \left\lbrace S_t \right\rbrace_{t >1} \subset \mgn$ converging to $C$, the sequence of maps 
$ \frac{1}{\log(t)} \cdot \mathcal{A}_R : S_t \rightarrow \mathbb{R}^m $ converges to the map $\pi_R : C \rightarrow \mathbb{R}^m$, in the sense of Definition \ref{GHconv}.
\end{mytheorem}

The notion of convergence of Definition  \ref{GHconv}
implies in particular that $\frac{1}{\log(t)} \cdot \mathcal{A}_R(S_t)$ converges towards $\pi_R(C)$ in Hausdorff distance on compact sets. 

Theorem \ref{approxtrop} will arise as a consequence of Theorem \ref{convdif} describing the limit of  imaginary normalised differentials  on families of algebraic curves converging in the sense of Definition \ref{tropconv}, see Section \ref{sec:convind}.
\medskip

\subsubsection{Phase-tropical convergence and Mikhalkin's approximation Theorem.}$\,$\medskip

Phase-tropical curves play a central role in the computation of planar Gromov-Witten invariants, see \cite{Mikh05}. They are complexified version of tropical curves and arise as Hausdorff limits of families of curves in $(\bC^\star)^m$ after the re-parametrisation
\[
\begin{array}{rcl}
H_t \; \;  : \;\; (\bC^\star)^m & \rightarrow & (\bC^\star)^m \\
(z_1,\dots,z_m) & \mapsto & \left( \vert z_1 \vert^{\frac{1}{\log(t)}} \frac{z_1}{\vert z_1 \vert} ,\dots, \vert z_m \vert^{\frac{1}{\log(t)}} \frac{z_m}{\vert z_m \vert} \right)
\end{array}.
\]
Classically, a family of curves $\left\lbrace \C_t \right\rbrace_{t>1} \subset (\mathbb{C}^\star)^m$ is said to converge to a phase-tropical curve $V \subset (\mathbb{C}^\star)^m$ if 
\begin{equation}\label{eq:convpt}
\lim_{t\rightarrow \infty} \;  H_t (\C_t) = V
\end{equation}
with respect to the Hausdorff distance on compact sets. 
In general, phase-tropical curves in $(\bC^\star)^m$ can be very singular objects. For this reason, we will restrict our attention to simple phase-tropical curves in this paper.

Simple phase-tropical curves can be considered abstractly. An abstract simple phase-tropical curve $V$ can be constructed by gluing (phase-tropical) pairs of pants together and comes with a map $\A_V:V\rightarrow C$ onto a simple tropical curve $C=(G,\ell)$. We will see in Section \ref{sec:trop} that any such object can also be described in terms of Fenchel-Nielsen coordinates and be therefore denoted by $V:=V(G,\ell,\theta)$.

\begin{definition}[\textbf{Abstract phase-tropical convergence}]\label{comptropconv}$\,$\\
Let $V:=V(G,\ell,\theta)$ be a simple phase-tropical curve such that $C=(G,\ell)$ is a simple tropical curve of genus $g$ with $n$ leaves. 
A family of Riemann surfaces $ \left\lbrace S_t \right\rbrace_{t >1} \subset \mgn$ \df{converges to} $V$ if $ S_t := S(G,\ell_t, \theta_t)$ for some functions $(\ell_t,\theta_t):E(G) \rightarrow \R_{> 0}\times S^1$ and if for any edge $e\in E(G)$, we have 
$$ \displaystyle \ell_t(e) \;\;  \underset{t\rightarrow \;\;  \infty}{\sim} \frac{2\pi^2}{\ell(e) \log(t)}  \; \;\; \text{ and } \;\;\;  \lim_{t\rightarrow \infty} \theta_t(e)=\theta(e).$$
\end{definition}

The above notion of convergence leads to a refined version of Theorem \ref{convdif}. For any $m \times n$ matrix of residues $R$, we can consider the \df{period matrix} $ \cP_{R,S}$ of the $m$ i.n.d.  on $S \in \mgn$ induced by $R$, where the periods are computed against a fixed basis of $H_1\big( S, \mathbb{Z} \big)$. 

\begin{mytheorem}
For any sequence $\left\lbrace S_t \right\rbrace_{t>1} \subset \mgn $ converging  to a simple phase-tropical curve $V$ and any $m \times n$ matrix of residues $R$, the sequence of period matrices $\mathcal{P}_{R,S_t}$ converges to a matrix $\mathcal{P}_{R,V}$ depending only on $R$ and $V$.
\end{mytheorem}

The above theorem will be stated more precisely as Theorem \ref{thm:convperiod} in Section \ref{sec:cop}. Now, for an abstract phase-tropical curve $V$ with underlying tropical curve $C$, any phase-topical morphism $\phi:V\rightarrow (\bC^\star)^m$ is induced from a tropical morphism $\pi:C\rightarrow\R^m$ induced in turn from an $m\times n$ matrix of residues $R$, see Propositions \ref{prop:phasetropdetermined} and \ref{propharmmorph}. By Theorem \ref{approxtrop}, the family of maps $\frac{1}{\log(t)}\cdot \A_R:S_t\rightarrow\R^m$ converges to $\pi:C\rightarrow\R^m$ for any sequence $\left\lbrace S_t \right\rbrace_{t>1}$ converging to $C$. In order to approximate the phase-tropical morphism $\phi$, we need to guarantee that the maps $\A_R:S_t\rightarrow\R^m$ factorise through algebraic maps $\phi_t:S_t\rightarrow(\bC^\star)^m$, or equivalently that the period matrices $\cP_{R,S_t}$ have coefficient in $2\pi i\bZ$. This will follow from Theorem \ref{thm:convperiod}, leading to a new proof of Mikhalkin's approximation Theorem, stated as Theorem \ref{thmMik} in Section \ref{Mik}.

\medskip

\subsection{Techniques and perspectives}

The main ingredient behind Theorems \ref{convdif} and \ref{thm:convperiod} is a rather elementary study of the behaviour of meromorphic differentials on families of annuli and pairs of pants respectively, see Sections \ref{secconfinv} and \ref{sec:cop}. The latter study provides us with a good understanding of the behaviour of i.n.d. on families of Riemann surfaces. In particular, we could design the Definitions \ref{tropconv} and \ref{comptropconv} of tropical and phase-tropical convergence for which the desired approximations theorems hold, namely Theorems \ref{approxtrop} and \ref{thmMik}. 

In this paper, we restrict our attention to simple (phase-) tropical curves, mostly for the sake of clarity and concision. It is of primary interest to extend the present results to general tropical curves. The most interesting application is probably related to compactifications of the moduli spaces $\mgn$ with tropical curves. Such compactifications are promising candidates to provide systematic correspondence theorems between algebro-geometric problems and their tropical counterpart.

From this perspective, Theorems \ref{approxtrop} and \ref{thmMik} justify the relevance of Definitions \ref{tropconv} and \ref{comptropconv}. Yet an other justification is given in Remark \ref{rem:jacobian} were we describe the relation between the Jacobians of a family of algebraic curves converging tropically and the Jacobian of the tropical limit, as defined in \cite{MZ}. The latter observation suggests a compactification of the Torelli map by its tropical analogue, as constructed in \cite{BMV}.

To conclude, let us mention that similar compactifications of moduli spaces have been investigated in \cite{O} and \cite{O2}. However, the present approach is quite different, as explained in Remark \ref{rem:odaka}. 

\tableofcontents

\section{Prerequisites}\label{sec:prereq}

\subsection{A twisted Hodge bundle}\label{sec:thb}

In order to study the convergence of imaginary normalised differential over families of Riemann surfaces, we will need an appropriate ambient space.

Recall that the Hodge bundle $\Lambda_{g,n} \rightarrow \bmgn$ is the rank $g$ vector 
bundle whose fiber at a point $(S;p_1,\dots,p_n)$ is the space of holomorphic sections of the dualising sheaf over $S$, see 
\cite{ELSV}. Geometrically, the fiber of $\Lambda_{g,n}$ over a smooth curve $(S;p_1,\dots,p_n) \in \mgn$ is 
the vector space of holomorphic differentials on $S$. For a singular curve $(S;p_1,\dots,p_n)  \in 
\bmgn$, the fiber of $\Lambda_{g,n}$ over $S$ is the vector 
space of meromorphic differentials whose pullback to the normalisation of $S$ have at 
most simple poles at the preimages of the nodes and such that the residues at the two preimages of any node of $S$ are opposite to each other. 
%

\begin{definition}
Define the \df{twisted Hodge bundle} $\thb_{g,n} \rightarrow \bmgn$ to be the vector bundle of rank $g+n-1$ obtained by tensoring $\Lambda_{g,n}$ with $\mathcal{O}_S (p_1+\dots+p_n)$.
A point in the fiber of $\thb_{g,n}$ over $S$ will be called a \df{generalised meromorphic differential} on $S$.
\end{definition}

Compared to meromorphic differentials, generalised meromorphic differential on $S$ are allowed to have simple poles at the punctures of $S$.

Let us now introduce coordinates on the fiber of $\thb_{g,n}$ over a fixed curve $S \in \mgn$.
Choose a maximal collection of points $q_1,\dots,q_k$ among the nodes of $S$ such that $S\setminus \big(\cup_j q_j \big)$ is connected. Choose now a maximal collection of simple closed curves $\gamma_1,\dots,\gamma_{g'}$ in $S$ such that the complement of these curves in $S\setminus \big(\cup_j q_j \big)$ is connected (observe that $g=g'+k$). Then, any generalised meromorphic differential $\omega$ on $S$ is uniquely determined by the vector 
\begin{equation}\label{eq:vector}
\textstyle \big( \int_{\gamma_1} \omega,\dots,\int_{\gamma_{g'}}\omega, \Res_{p_1}\omega, \dots,\Res_{p_{n-1}}\omega, \Res_{q_1}\omega, \dots,\Res_{q_{k}}\omega\big) \in \bC^{g+n-1}.
\end{equation}
Observe that we need to choose a branch of the node $q_j$ in order to define $\Res_{q_j}$.

We can also consider similar coordinates along families of algebraic curves.
For any continuous family $\left\lbrace S_t \right\rbrace_{t>1}\subset \mgn$ converging to $S$, each node $q_j$ defines an isotopy class of simple closed curve $\gamma'_{j,t}\subset S_t$, the so-called vanishing cycle associated to $q_j$. Now, a family $\left\lbrace (S_t,\omega_t) \right\rbrace_{t>1}\subset \thb_{g,n}$ converges to the point $(S,\omega)$ if and only if the vector 
\[
\textstyle \big( \int_{\gamma_1} \omega_t,\dots,\int_{\gamma_{g'}}\omega_t, \Res_{p_1}\omega_t, \dots,\Res_{p_{n-1}}\omega_t, \int_{\gamma'_1}\omega_t, \dots,\int_{\gamma'_k}\omega_t\big)
\]
converges to the vector \eqref{eq:vector}. Observe that $\gamma'_{j,t}$ has to be oriented coherently with the choice we made to define $\Res_{q_j}$.

\subsection{Imaginary normalised differentials and harmonic amoebas}\label{sec:indha}

In this section, we briefly review the material of \cite{Kri} that is of interest to us. Unless specified otherwise, the proofs of the statements to follow can be found there.

\begin{definition}
An \df{imaginary normalised differential} (\df{i.n.d.} for short) $\omega$ on a curve $S \in \bmgn$ is a generalised meromorphic differential having (at worst) simple poles at the $n$ punctures of $S$ and such that for any simple closed curve $ \gamma \subset S$, we have
\[ \Re \left( \int_{\gamma} \omega \right) = 0. \]
\end{definition}

\begin{theorem}\label{thmRiem}
Let $S \in \bmgn$ with punctures $p_1,\dots, p_n$ and nodes $q_1,\dots, q_k$. Fix any collection $\left\lbrace r_1,\dots, r_n\right\rbrace \subset \R$ such that $\, \sum_j r_j =0$ and associate a number $ r'_j \in \R$ to one of the two branches of the node $q_j$. Then, there exists a unique i.n.d. $\omega$ on $S$ such that $\, \Res_{p_j} \omega = r_j \,$ for any $1 \leqslant j \leqslant n$ and $\, \Res_{q_j} \omega = r'_j \,$ for any $1 \leqslant j \leqslant k$ where the latter residue is computed at the distinguished branch of $q_j$.
\end{theorem}

\begin{proof}
For $S \in \mgn$, the proof is given in \cite[Theorem 5.3]{Lang82}. In the general case, normalise the curve $S$, associate the residue $-r'_j$ to the remaining branch of the node $q_j$ and apply \cite[Theorem 5.3]{Lang82} to each connected component of the normalisation.
\end{proof}

\begin{definition}\label{defAr}
A \df{collection of residues} of dimension $m$ on a Riemann surface $S \in \mgn$ is a real $n \times m$ matrix $R:= \big( r_{k,j} \big)_{k,j}$ such that $\sum_k r_{k,j} =0$ for any $1 \leqslant j \leqslant m$. Denote by $$\omega_R^S :=(\omega_{R,1}^{S}, \dots, \omega_{R,m}^{S})$$ the vector of i.n.d. on $S$ induced by the $m$ rows of $R$, see Theorem \ref{thmRiem}. In practice, we will simply denote the latter vector by $\omega_R=(\omega_{R,1}, \dots, \omega_{R,m})$ when no confusion is possible. Finally, define the \df{harmonic amoeba map} 
\[
\begin{array}{rcl}
\mathcal{A}_R \;\;  : \;\;   S & \rightarrow & \mathbb{R}^m \\
s & \mapsto & \Big( \Re\big( \int_p^s \omega_{R,1}\big), \dots, \Re\big( \int_p^s \omega_{R,m}\big)\Big)
\end{array}
\]
%
up to the choice of an initial point $p \in S$. We call the set $\mathcal{A}_R  (S) \subset \mathbb{R}^m$ the \df{harmonic amoeba} of $S$ with respect to $R$.
\end{definition}

\begin{remark}
According to Theorem \ref{thmRiem}, the space of harmonic amoebas maps on any given Riemann surface $S\in \mgn$ to $\mathbb{R}^m$ is a real vector space of dimension $m(n-1)$ given by the coefficients of the collection of residues.
\end{remark}

The terminology of Definition \ref{defAr} is motivated by the fact that the coordinates of $\mathcal{A}_R$ are harmonic functions on $S$ and that the definition of harmonic amoebas generalises the one of \cite{GKZ}, as we recall below. 

The next proposition illustrate the fact that the main properties of amoebas survive in the context of harmonic amoebas, see \cite[Proposition 2.7]{Kri}.

\begin{proposition}\label{convamoeba}
Let $S \in \mgn$ and $R$ be a collection of residues of dimension 2. Then, the harmonic amoeba $\mathcal{A}_R  (S)$ is a closed subset of $\R^2$ with finite area and the connected components of $\mathbb{R}^2 \setminus \mathcal{A}_R  (S)$ are convex.
\end{proposition}

It is also shown in \cite{Kri} that harmonic amoebas possess a logarithmic Gauss map, a Ronkin function and extend many classical properties of amoebas. In particular, we can carry the definition of spine of planar amoebas as introduced in \cite{PR} to the case of harmonic amoebas. The latter construction suggests to consider a more general class of affine tropical curve with non rational slopes that we introduce in section \ref{secharmtrop} under the name of  \df{harmonic tropical curves}.

Observe that for any collection of residues $R$ of dimension $m\geqslant 2$, the image of any small disc $D\subset (S\cup p_j)$ around the puncture $p_j$ is mapped by $\A_R$ inside the $\varepsilon$-neighbourhood of a half-ray $r\subset \R^m$ and that $\varepsilon$ can be made arbitrarily small by shrinking $D$. Indeed, if we consider a holomorphic coordinate $t$ on $D$ centred at $p_j$ and an initial point $p\in D$, we have that 
\[
\begin{array}{rl}
\A_R(s) & = \; \Big( \Re\big( \int_p^s \frac{r_1}{t}+ O(1) dt\big), \dots, \Re\big( \int_p^s \frac{r_m}{t}+ O(1) dt\big)\Big)\\ & \\
&  = \; \big( r_1 \log\vert s \vert  + O(1), \dots, r_m \log \vert s \vert  + O(1)\big).
\end{array}
\]
The vector $(r_1,\dots,r_m)^\top$ is the $j^{th}$ column of $R$, each $O(1)$ is meant in a neighbourhood of $0$ and can be made as small as desired by shrinking $D$. Following the terminology of \cite{GKZ} and \cite{Kri}, we refer to $\A_R( D\setminus p_j)$ as a \df{tentacle} of the harmonic amoeba $\A_R(S)\subset \R^m$ and to $(r_1,\dots,r_m)^\top$ as the \df{slope} of the tentacle. 

For $S \in \mgn$, a collection of residues $R$ of dimension $m$ and the integer $k:=2g+n-1$, we can define the (normalised) \df{period matrix} $ \mathcal{P}_{R,S} \in M_{k\times m}(\bC)$ relative to a basis $\gamma_1,\dots, \gamma_{k}$ of $H_1(S,\mathbb{Z})$ as follows
\[  \big( \mathcal{P}_{R,S} \big)_{j,k} := \frac{1}{2\pi i}  \int_{\gamma_j} \omega_{R,k} \; . \]
The period matrix relative to a different basis $\gamma':= A \cdot \gamma$  is given by $A \cdot \mathcal{P}_{R,S}$ where $A \in \SL_{k}(\mathbb{Z})$.

\begin{definition}\label{defpermat}
For $S \in \mgn$ and $R$ be a collection of residues of dimension $m$, the \df{period matrix} of $S$ with respect to $R$ is the equivalence class of the matrix $\mathcal{P}_{R,S}$ relative to any basis of $H_1(S,\mathbb{Z})$ in the quotient of $M_{k\times m}(\bC)$ by $\SL_{k}(\mathbb{Z})$ acting by multiplication to the left. By a slight abuse of notation, we denote the latter class by $\mathcal{P}_{R,S}$. The class $ \mathcal{P}_{R,S} $ is an \df{integer period matrix} if one of its representative (and then all of them) has coefficients in $\mathbb{Z}$. If $ \mathcal{P}_{R,S} $ is an integer period matrix, define the holomorphic map 
\[
\df{\begin{array}{rcl}
\mathcal{\phi}_{R}  \; : \; S & \rightarrow & (\bC^\star)^m\\
z & \mapsto & \left( e^{\int^z_{p} \omega_{R,1}},\dots,  e^{ \int^z_{p} \omega_{R,m}}  \right)
\end{array}}
\]
up to the choice of an initial point $p\in S$
\end{definition}

The next statement follows from Remark 3.2 and Theorem 2.6 in \cite{GK}.

\begin{theorem}\label{thmcodim}
Let $S \in \mgn$ and $R$ be a collection of residues of dimension$ \; 1$. The level set
\[\left\lbrace S' \in \mgn \: \big| \: \mathcal{P}_{R,S} = \mathcal{P}_{R,S'} \right\rbrace\]
is a smooth analytic subvariety of $\mgn$ of real codimension $2g$.
\end{theorem}\medskip

\subsection{Hyperbolic surfaces}\label{sec:hyp}

In this section, we recall some basic facts about geometry of hyperbolic surfaces. We refer to \cite{B} and \cite{Hub} for more details.

\begin{definition}
A \df{pair of pants} $Y$ is a complete Riemannian surface with constant curvature $-1$ and  geodesic boundary that is conformally equivalent to  $ \bC P^1 \setminus \big( E_1 \cup E_2 \cup E_3 \big) $ where $E_j$ is either a point or an open disc and such that the $E_j$'s are pairwise disjoint.
\end{definition}

Recall  that the upper half-space $ \mathbb{H} := \left\lbrace x+iy \in \mathbb{C} \; \big| \;  y>0 \right\rbrace $ equipped with the Riemannian metric 
\[g:= \dfrac{(dx^2 +dy^2)}{y^2}\]
is a complete metric space with constant curvature $-1$, see \cite[Ch.1, \S 1]{B}.

\begin{theorem}[Theorem 2.4.2 in \cite{B}]
For any $\alpha, \beta, \gamma >0$, there exists a hyperbolic right-angled hexagon $H(\alpha, \beta, \gamma) \subset \mathbb{H}$ unique up to isometry such that its sides have consecutive length $a, \alpha, b, \beta, c, \gamma$. The lengths $a,b,c$ are determined by $\alpha, \beta, \gamma$.
\end{theorem}

We can define generalised hyperbolic right-angled hexagon by letting any of the parameters $\alpha, \beta, \gamma$ go to $0$. For any $(\alpha, \beta, \gamma)\in (\R_{\geqslant 0})^3$, we obtain a Riemannian surface $H(\alpha, \beta, \gamma) \subset \mathbb{H}$ unique up to isometry.

From now on, by a boundary component of a pair of pants $Y$ we mean a boundary geodesic as well as a puncture of $Y$. The theorem below is a combination of Proposition 3.1.5, Theorem 3.1.7 and Lemma 4.4.1 in \cite{B}.

\begin{theorem}\label{thm:buser}
For any pair of pants $Y$, there exists a unique orientation reversing isometry $\sigma$ fixing globally each boundary component. The quotient $Y/\sigma$ is isometric to a generalised hyperbolic right-angled hexagon. In particular, the surface $Y$ is determined up to isometry by the length $(\alpha, \beta, \gamma)\in (\R_{\geqslant 0})^3$ of its $3$ boundary components.
\end{theorem}

The isometry $\sigma$ in the above theorem is an involution. By analogy with complex conjugation, we will denote by $\R Y \subset Y$ the \df{fixed locus of $\sigma$}.

Now let us recall the construction of Riemann surfaces using Fenchel-Nielsen coordinates. Below, we use $l(\gamma)$ to denote the length of a closed geodesic in a hyperbolic surface and $d(\, . \, ,\, . \,)$ denotes the distance induced by the hyperbolic metric.

\begin{definition}\label{def:cubgraph}
A \df{cubic graph} $G$ is a connected graph with only trivalent and univalent vertices. An edge adjacent to a univalent vertex is called a \df{leaf}. A \df{ribbon structure} $\mathscr{R}$ on $G$ is the data for every vertex $v$ of $G$ of a cyclical ordering of the edges adjacent to $v$. We denote by $V(G)$ the \df{set of vertices}, by $L(G) $ the \df{set of leaves} and by $E(G)$ the \df{set of non-leaf edges} of $G$. We also denote $LE(G):= L(G) \cup E(G)$.
\end{definition}

Let us fix a cubic graph $G$ with $n$ leaves and genus $g$, that is $g:=b_0(G)$, together with a fixed ribbon structure $\scrR$. By Euler characteristic computation, we have that $G$ has $(3g-3+n)$ edges and $(2g-2+n)$ vertices of valence $3$. We assume that $2g+n>2$.

Given a \df{length function} $\ell \, : \, E(G)  \rightarrow  \mathbb{R}_{> 0} $  and a \df{twist function} $\theta \, : \, E(G)  \rightarrow  S^1$, we can construct a Riemann surface $S(G,\ell,\theta) \in \mgn$ as follows. To each vertex $v\in V(G)$, associate a pair of pants $Y_v$ together with a correspondence between the $3$ edges of $G$ adjacent to $v$ and the $3$ boundary components of $Y_v$ such that :

-- a boundary component of $Y_v$ is a puncture if and only if the corresponding edge of $G$ is a leaf,

-- the boundary geodesic of $Y_v$ corresponding to the edge $e\in E(G)$ has length $\ell(e)$.\\
Now, we will identify each boundary geodesic $\gamma\subset Y_v$ to $S^1$ with the help of the ribbon structure $\scrR$. Any pointed, oriented topological circle equipped with an inner metric of length $2\pi$ is uniquely isomorphic to $S^1$, pointed at $1$, oriented counter-clockwise and equipped with the Euclidean metric. Thus, 
we simply need to choose an ``origin" and an orientation on $\gamma$, equipped with the rescaled hyperbolic distance $\frac{2\pi}{l(\gamma)}d(\,.\,,\, .\,)$, in order to identify $\gamma$ with $S^1$. We choose the orientation whose direct normal vector is pointing inside $Y_v$. For the origin, observe that the involution $\sigma$ on $Y_v$ has two fixed points on $\gamma$. Exactly one of these two points belong to a component of $\R Y_v$ intersecting with the boundary component of $Y_v$ preceding $\gamma$ according to the ordering induced by $\scrR$. We choose the latter point as the origin of $\gamma$.

It remains to glue the pairs of pants $Y_v$ together. For two vertices $v_1, v_2 \in V(G)$ connected by an edge $e\in E(G)$, glue the two geodesics  of $ Y_{v_1} $ and $ Y_{v_2} $ corresponding to $e$ with the isometric involution 
\[
\begin{array}{rcl}
 S^1 & \rightarrow &  S^1\\
z & \mapsto & -\overline{\theta(e)z}.
\end{array}\]
The result of the above construction is a Riemann surface $S(G,\ell,\theta) \in \mgn$. The functions $(\ell,\theta)$ are the \df{Fenchel-Nielsen coordinates} of the Riemann surface $S(G,\ell,\theta)$, relative to $G$ and $\scrR$. These coordinates allow to define the map
\[
\begin{array}{rcl}
\fn_G \: : \: (\bC^\star)^{3g-3+n} & \rightarrow & \mgn \\
(\ell,\theta) & \mapsto & S(G,\ell,\theta)
\end{array}
\]

\begin{theorem}\label{thmhyperb}
The map $\fn_{G}$ is a real-analytic surjective map.
\end{theorem}

\begin{proof}
Consider the Teicmüller space $\mathcal{T}_{g,n}$ together with its canonical complex structure so that the universal map $\mathcal{T}_{g,n}  \rightarrow \mgn$ is holomorphic. The Teichmüller space admits Fenchel-Nielsen coordinates $(\ell, \Theta)\in (\R_{>0})^{3g-3+n}\times \R^{3g-3+n}$ (relative to any choice of $G$ and $\scrR$) such that the length function $\ell$ is the same as the one we just defined and the twist function $\theta$ is the projection of $\Theta$ in $(\R/\bZ)^{3g-3+n}$, see \cite[Ch.3, \S 3]{B}. It follows that the universal map is the composition of $(\ell,\Theta)\mapsto(\ell,\theta)$ with $\fn_G$. In particular, $\fn_G$ is surjective since the universal map is. Also, we deduce that $\fn_G$ is real-analytic if and only the universal map is real analytic in Fenchel-Nielsen coordinates. According to \cite[Theorem 3.12]{Wol}, the Fenchel-Nielsen coordinates $(\ell, \Theta)$ are real-analytic on $\mathcal{T}_{g,n}$. In particular, the universal map $\mathcal{T}_{g,n}  \rightarrow \mgn$ is again analytic in the latter coordinates. The result follows.
%
\end{proof}

To conclude this section, we recall some classical facts around the collar Theorem, see \cite[Ch.4]{B} for more details. 

\begin{definition}
Let $Y$ be a pair of pants. For any boundary geodesic $\gamma$ of $Y$, the \df{half-collar} associated to $\gamma$ is the tubular neighbourhood
\[\hc_\gamma := \left\lbrace  z \in Y \; \big| \; d(z, \gamma) \leqslant w \big(l(\gamma)\big) \right\rbrace\]
where $ w(l) := \arcsinh \big(\frac{1}{\sinh( l/2)}\big)$. For any puncture $p$ of $Y$, the \df{cusp} $K_p$ associated to $p$ is the neighbourhood of $p$ in $Y$ isometric to the space $\left]- \infty, \log(2) \right] \times \mathbb{R}/\mathbb{Z}$ with  coordinates $(\rho,y)$ and metric $d\rho^2 + e^{2 \rho}dy^2$. 
For any Riemann surface $S \in \mgn$, and any simple closed geodesic $\gamma \subset S$, the \df{collar} associated to $\gamma$ is the tubular neighbourhood
\[K_\gamma := \left\lbrace  z \in S \; \big| \; d(z, \gamma) \leqslant w\big(l(\gamma)\big) \right\rbrace.
\]
\end{definition}

According to \cite[Equation (3.3.6)]{B}, the half-collar $\hc_\gamma$ is isometric to $\left[0, w(l) \right] \times \mathbb{R}/\mathbb{Z}$ equipped with the metric $d\rho^2 + l^2 \cosh(\rho)^2 dy^2$.

\begin{lemma}\label{propcolcusp}
The half-collar to a geodesic $\gamma$ converges to a cusp (in Gromov-Hausdorff distance) when $l(\gamma)$ tends to $0$.
\end{lemma}

\begin{proof}
To see this, it suffices to compare the limit of the length of the boundary component of $\hc_{\gamma}$ different from $\gamma$ with the length of the boundary of $K_p$. The latter lengths are easily computed from the metric descriptions given above. The length is $2$ for the cusp and $l\cosh(w(l))$ for the half-collar. Using the formula $\arcsinh(x)=\ln(x+\sqrt{x^2+1})$ and the fact that $\lim_{l \rightarrow 0} w(l)= +\infty$, we compute that
\[
\begin{array}{rl}
l^2\cosh(w(l))^2 & = \;\,  l^2\Big(\frac{e^{w(l)}+e^{-w(l)}}{2} \Big)^2 \;\,  \underset{l \rightarrow 0}{\sim}\;\,   l^2\Big(\frac{e^{w(l)}}{2} \Big)^2 \;\,  \sim \;\,  l^2\left(\frac{\frac{1}{\sinh(l/2)}+\sqrt{\frac{1}{\sinh(l/2)^2}+1}}{2} \right)^2 \\ \\
& \sim \;\,  l^2\Big(\frac{1}{\sinh(l/2)} \Big)^2 \;\,  \sim \;\,  l^2 \big(\frac{2}{l}\big)^2 \;\, = \;\,  4 \, .
\end{array}
\]
\end{proof}

\begin{proposition}[Proposition 4.4.4 in \cite{B}]
Let $Y$ be a pair of pants. The half-collars and cusps of $Y$ are pairwise disjoint.
\end{proposition}
\medskip

\subsection{Conformal invariants on holomorphic annuli}\label{secconfinv}

Unlike holomorphic discs, holomorphic annuli have a non trivial modulus. For example, the annuli $A_{r,R} = \left\lbrace z \in \mathbb{C} \: \big| \: r \leqslant \vert z \vert \leqslant R \right\rbrace$ are determined up to bi-holomorphism by the single conformal invariant $\frac{1}{2\pi}\log \left( \frac{R}{r}\right)$, see \cite[Ch.14]{Rud}. This conformal invariant is related to so-called extremal lengths on complex domains, see e.g. \cite{Ahl}. Since any holomorphic annulus is bi-holomorphic to $A_{1,R}$ for some $R>0$, holomorphic annuli have a single modulus, see either \cite[Ch.5, \S 1, Theorem 10]{Ahl79} or \cite[Proposition 3.2.1]{Hub}. An other useful model for holomorphic annuli is given by 
\[A_m = \left\lbrace z \in \mathbb{C} \: \big| \: 0 \leqslant \Re( z ) \leqslant m \right\rbrace \big/ i \mathbb{Z},\]
for any $m>0$.
The map $ z \mapsto e^{2\pi z} $ identifies $A_m$ to $A_{1,e^{2\pi m}}$ so that $A_m$ has modulus $m$. As a consequence, every holomorphic annulus is bi-holomorphic to  $A_m$ for a unique $m>0$.

In the sequel, we will need to compare the above conformal models with the collar model $K_l$ coming from hyperbolic geometry.

\begin{lemma}\label{lemmaconfcollar}
Define $m(l):= \frac{2}{l}\arccos\big(\frac{1}{\cosh(w(l))}\big)$. Then, the collar $K_l$ is bi-holomorphic to the annulus $A_{m(l)}$. In particular, we have
\[ m(l) \; \underset{l \rightarrow 0}{\sim} \;  \pi/l \,  .\]
\end{lemma}

\begin{proof}
Recall that the collar $K_l$ is isometric to $\left[- w(l), w(l) \right] \times \mathbb{R}/\mathbb{Z}$ equipped with the metric $d\rho^2 + l^2 \cosh(\rho)^2 dy^2$. Let us look for a change of variable $x \mapsto \rho(x)$ so that the pullback of the latter metric  to the $(x,y)$-space is conformally equivalent to the Euclidean metric. This amounts to solving the differential equation
\[\big(\rho'(x)\big)^2 = \big(l  \cosh (\rho(x))\big)^2.\]
A solution is given by the odd function $\rho(x)$ whose restriction to $\R_{\geqslant 0}$ is given by $\rho(x) = \arccosh\big(\frac{1}{\cos(lx)}\big)$.
It follows that $K_l$ is isometric to $\left[ - \rho^{-1}\big(w(l)\big), \rho^{-1}\big(w(l)\big) \right] \times \mathbb{R}/\mathbb{Z}$ with the metric 
$\big(\frac{l}{\cos(lx)}\big)^2 (dx^2 + dy^2)$. The latter is the hyperbolic metric on the annulus $A_{2\rho^{-1}(w(l))}$. We deduce that $$m(l) \; = \;  2\, \rho^{-1}\big(w(l)\big)\; = \;  \frac{2}{l}\,\arccos\left(\frac{1}{\cosh( w(l))}\right).$$
For the asymptotic of $m(l)$ at $0$, observe  that  $$\displaystyle \lim_{l\rightarrow 0} \,  w(l)  =   + \infty \; \; \text{ and } \; \;  \displaystyle \lim_{x\rightarrow +\infty} \, \arccos\left(\frac{1}{\cosh(x)}\right) = \frac{\pi}{2}.$$
\end{proof}

In the sequel, we will also need to consider families of holomorphic differentials 
$\left\lbrace \omega_m \right\rbrace_{m>0}$ where $\omega_m$ is defined on $A_m$. In order to deal with such families, it is useful to introduce the model
\[C_\lambda = \left\lbrace (z,w) \in \mathbb{C}^2 \, \big|  \; \vert z \vert \leqslant 1,  \; \vert w \vert \leqslant  1, \: zw= \lambda \right\rbrace\]
for $0 \leqslant \lambda <1$. The parametrisation $z \mapsto (e^{-2\pi z},\lambda e^{2\pi z})$ identifies $C_\lambda$ with $A_{-\log(\lambda)/(2\pi)}$. In particular, the family $\big\lbrace C_\lambda\big\rbrace_{0<\lambda<1}$ is a universal family of annuli. For $\lambda=0$, the set $C_\lambda$ is the union of the 2 unit discs in the $z$- and the $w$-coordinate axes of $\bC^2$. 

\begin{definition}\label{deflocconv}
A family  $\left\lbrace \omega_m \right\rbrace_{m>0}$ of holomorphic differentials defined respectively on $A_m$ is \df{continuous} if the differential induced on the total space $\cup_{0<\lambda<1} C_\lambda \subset \bC^2$ is continuous. A continuous family $\left\lbrace \omega_m \right\rbrace_{m>0}$ \df{converges} at $+\infty$ if the limit in the $C_\lambda$-model is a meromorphic differential on each of the two  irreducible components of $C_0$.
\end{definition}

Observe that the limiting differential is holomorphic everywhere on $C_0$ except at the origin where the residue on the two discs are opposite to each other. Suppose now that we are given a continuous family of holomorphic differentials 
$\left\lbrace \omega_m \right\rbrace_{m>0}$ as above and normalised in the following way
\begin{eqnarray}\label{normdiff}
\dfrac{1}{2 \pi i} \int_\gamma \omega_m = 1
\end{eqnarray}
where $\gamma(t)=m/2+it$, $0\leqslant t \leqslant 1$.
We would like to compute the asymptotic of the integral of $\omega_m$ along the width of $A_m$ while $m$ tends to $+\infty$. To that aim, let us introduce the following terminology. 

\begin{definition}
An admissible sequence of path is a sequence of continuous, piecewise $C^1$ maps $\rho_m :  [0,m] \rightarrow A_m$ such that

-- $\Re \big( \rho_m (0) \big) =0$ and $\Re \big( \rho_m (m) \big) =m$,

-- there exists $M>0 $ such that $ \left| \int_0^m \frac{d}{dt}\Im\big( \rho_m(t) \big) dt \right| <M$ for any $m$.
\end{definition}

The latter condition is equivalent to saying that the parametrised curve obtained by composition of $\rho_m$ with the projection $A_m\rightarrow \R/\bZ$ has length at most $M$. 

\begin{lemma}\label{lemmaconfdiff}
Let $\left\lbrace \omega_m \right\rbrace_{m>0}$ be a converging family of holomorphic differentials normalised as in \eqref{normdiff}. For any admissible sequence of paths $\left\lbrace \rho_m \right\rbrace_{m>0}$, we have
\[\int_{\rho_m} \omega_m \; \underset{m \rightarrow \infty}{\sim} \; 2 \pi m.\]
\end{lemma}

\begin{proof}
First, consider the family $\omega_m^0 = 2 \pi dz$ on $A_m$, obviously normalised as in \eqref{normdiff}. Recall that the map $z \mapsto (e^{-2\pi z},\lambda e^{2\pi z})$ allows to identify $C_{e^{-2\pi m}}$ with $A_{m}$. It follows then from a simple computation that the differential  $\omega_m^0$ is given on each $C_{e^{-2\pi m}}$ by the restriction of the differential
\[\left(-\frac{1}{2z}+w\cdot h(z,w)\right) dz + \left(\frac{1}{2w}+z\cdot h(z,w)\right) dw\]
on $\bC^2$, where $h(z,w)$ is any meromorphic function defined on $(\bC^\star)^2$. In particular, we see that $\omega_m^0$ converges to the meromorphic differential on $C_0$ that restricts respectively to $-dz/z$ and $dw/w$ on the $z-$ and $w-$unit discs by taking $h(z,w)=\pm (2zw)^{-1}$. Consider now the admissible sequence of paths $\rho_m^0 (t)=t$. Then, we have 
\[\int_{\rho_m^0} \omega_m^0 = 2 \pi m.\]
Now, for any converging family $\left\lbrace \omega_m \right\rbrace_{m>0}$ normalised as in \eqref{normdiff}, the sequence $\left\lbrace \omega_m - \omega_m^0 \right\rbrace_{m>0}$ converges to a differential on $C_0$ whose restriction to any component is holomorphic. It follows that $\int_{\rho_m^0} ( \omega_m - \omega_m^0 ) $ converges and that 
\[\int_{\rho_m^0} \omega_m \; \underset{m \rightarrow \infty}{\sim} \;  2 \pi m.\]
Now for any admissible sequence of paths $\left\lbrace \rho_m \right\rbrace_{m>0}$, there exist 2 sequences of  locally injective paths $\left\lbrace \rho_{1,m} \right\rbrace_{m>0}$ and $\left\lbrace \rho_{2,m} \right\rbrace_{m>0}$, respectively contained in the boundaries $ \left\lbrace \Re(z)=0 \right\rbrace$ and $ \left\lbrace \Re(z)=m \right\rbrace$ of $A_m$, and such that $\rho_m$ is isotopic to $\rho_{2,m} \circ \rho_m^0 \circ \rho_{1,m}$.
Hence, we have that
\[\int_{\rho_m} \omega_m = \int_{\rho_{1,m}} \omega_m + \int_{\rho_m^0} \omega_m + \int_{\rho_{2,m}} \omega_m .\]
Since $\left\lbrace \rho_m \right\rbrace_{m>0}$ is admissible, the lengths of the paths $\rho_{1,m}$ and $\rho_{2,m}$ are bounded from above. Since $\left\lbrace \omega_m \right\rbrace_{m>0}$ converges, the integrals of $\omega_t$ along $\rho_{1,m}$ and $\rho_{2,m}$ are also bounded from above. As a consequence, we have
\[\int_{\rho_m} \omega_m \; \underset{m \rightarrow \infty}{\sim} \;  \int_{\rho_m^0} \omega_m \; \underset{m \rightarrow \infty}{\sim} \;  2 \pi m.\]
\end{proof}
%

\begin{proposition}\label{convdifloc}
Let $\left\lbrace l_t \right\rbrace_{t>1}$ be a continuous family of positive numbers such that 
\[ l_t  \; \underset{t \rightarrow \infty}{\sim}  \; \dfrac{2\pi^2}{ l \cdot \log(t)} \]
for some $l>0$ and let $\alpha$ be a positive number. Then, for any converging sequence $\left\lbrace \omega_t \right\rbrace_{t>1}$ of holomorphic differentials defined respectively on $A_{\alpha \cdot m(l_t)}$ and any admissible sequence of paths $\left\lbrace \rho_t \right\rbrace_{t>1}$, we have 
\[\dfrac{1}{\log(t)} \int_{\rho_t} \omega_t  \; \underset{t \rightarrow \infty}{\sim} \;   \alpha \, l \, \Lambda\]
where 
\[  \Lambda \; := \; \lim_{t \rightarrow \infty} \; \dfrac{1}{2\pi i} \int_\gamma \omega_t. \] 
\end{proposition}

\begin{proof} Denote $\Lambda_t := 1/(2\pi i) \int_\gamma \omega_t$. Observe that, since the family $\left\lbrace \omega_t \right\rbrace_{t>1}$ converges, the integral $\Lambda_t$ converges to the residue of the limit differential on one of the two components of $C_0$. Assume first that $\Lambda \neq 0$ and let $t_0>1$ be such that $\Lambda_t \neq 0$ for any $t>t_0$.
Applying lemma \ref{lemmaconfdiff} to the normalised family $\left\lbrace (1/\Lambda_t)\omega_t \right\rbrace_{t>t_0}$, and the equivalence given in lemma \ref{lemmaconfcollar}, we obtain that
\[ \frac{1}{\Lambda_t} \int_{\rho_t} \omega_t \;\; \underset{t \rightarrow \infty}{\sim} \;\;  2 \pi \, \alpha \, m(l_t) \;\;  \underset{t \rightarrow \infty}{\sim} \;\;  2\pi \,\alpha\,  \frac{\pi}{l_t} \;\; \underset{t \rightarrow \infty}{\sim} \;\; \alpha\, l\, \log(t).\]
If  $\Lambda = 0$, then the limit differential is holomorphic on the two components of $C_0$ and $\int_{\rho_t} \omega_t$ is bounded as a function in $t$. The result follows.
\end{proof}
\medskip

\subsection{Tropical curves}\label{sec:trop}

In this section, we give a brief account on (phase-) tropical curves and tropical 1-forms on them.

\subsubsection{Tropical curves and morphisms}\label{sectropcurvmorph}

\begin{definition}\label{deftropcurv}
A \df{tropical curve} $C$ is a graph with all univalent vertices removed and equipped with a complete inner metric. Two tropical curves are \df{isomorphic} if they are isometric to each other. A tropical curve is \df{simple} if its underlying graph is a cubic graph (see Definition \ref{def:cubgraph}).
\end{definition}

We carry the notation $V(C)$, $L(C)$, $E(C)$ and $LE(C)$ from Definition \ref{def:cubgraph}. We recall that in the present terminology, edges have finite length while leaves have infinite length.
For any vertex $v \in V(C)$ of a simple tropical curve $C$, denote by $T_v \subset C$ the tripod obtained as the union of $v$ with the $3$ open leaves/edges adjacent to $v$ and $\htr_v \subset T_v$ the tripod obtained by taking only the half of each edge adjacent to $v$.
\begin{remark}\label{rem:tropGL}
A tropical curve can be represented by a couple $C :=(G,\ell)$, where $G$ is the graph supporting $C$ and $\ell :  E(G) \rightarrow \R_{>0} $ is the length  function induced by the metric on $C$.
\end{remark}
\begin{definition}\label{tropmorph}
A \df{tropical morphism} $ \pi \, : \, C \rightarrow \mathbb{R}^m $ on a tropical curve $C$ is a proper continuous map subject to the following 

--\textbf{Integrality}: for any $ e \in LE( C )$, the map $ \pi_{\vert e} $ is integer affine linear with respect to the metric on $e$.

--\textbf{Balancing}: For any $ v \in V ( C )$, denote by $\vec{e}_1, \dots, \vec{e}_n$ the  outgoing unitary tangent vectors to the $n$ leaves/edges adjacent to $v$. Then
\[ \sum_{1\leqslant j \leqslant n} d\pi (\vec{e}_{j}) \; = \; 0. \]
A \df{tropical curve} $\Gamma \subset \mathbb{R}^m$ is the image of tropical morphism $ \pi  :  C \rightarrow \mathbb{R}^m$.
\end{definition}

%

\subsubsection{Tropical 1-forms}

\begin{definition}\label{def1form}
A (tropical) \df{1-form} $\sw$ on a tropical curve $C$ is a locally constant, real-valued 1-form on $C\setminus V(C)$ that satisfies the following balancing condition: 

-- for any $v \in V( C )$, and outgoing unitary tangent vectors $\vec{e}_1,\dots,\vec{e}_n$ on the $n$ leaves/edges adjacent to $v$, we have the balancing condition 
\[ \sum_{1\leqslant j \leqslant n} \sw (\vec{e}_j) \; =\;  0. \]
The vector space of 1-forms on $C$ is denoted by $\Omega( C )$. For any inward oriented leaf $\vec{e}$ of $C$, define the \df{residue} of $\sw$ at $e$ to be the number $\sw (\vec{e})$. An element $\sw \in \Omega( C )$ is \df{holomorphic} if all its residues are zero. The vector space of holomorphic 1-form on $C$ is denoted by $ \Omega_\mathcal{H}( C )$.
\end{definition}

A 1-form $\sw$ on $C$ corresponds locally to the datum $a \, dx$ on any $e\in LE(C)$, where $a \in \mathbb{R}$ and $x \, : \, e \rightarrow \mathbb{R}$ is an isometric coordinate. For any other isometric coordinate $y$, one has $\sw = \pm a \, dy$ depending whether $x \circ y^{-1}$ preserves orientation or not. Therefore, the $1$-form $\sw$ is equivalent to the data of an orientation on any $e\in LE(C)$ plus the corresponding real number $a$. In other words, the form $\sw$ is equivalent to the datum of an electric current on $C$ seen as an electrical circuit. The above balancing condition corresponds then to the conservation of energy at the node of the circuit.

\begin{definition}\label{pathloop}
Let $C=(G,\ell)$ be a tropical curve. A \df{path} $C$ is a continuous map $\rho \, : \, [0,1] \rightarrow LE(C)$ such that $\rho(0)$, $\rho(1) \in V(G)$ and such that the restriction of $\rho$ to $[0,1[$ is injective. In particular, a path can join a univalent vertex of $G$ to an other. 
A \df{loop} is a path $\rho$ such that $\rho(0)=\rho(1)$. Alternatively, we will also represent a path by the induced ordered collection $(\vec{e}_1,\dots,\vec{e}_k)$ of oriented leaf/edge. For a path $\rho=(\vec{e}_1,\dots,\vec{e}_k)$ and any 1-form $\sw$ on $C$, we define the integral
\[\int_\rho \sw \; := \;  \sum_{1\leqslant j \leqslant k} \ell(e_j)\, \sw(\vec{e}_j). \]
A 1-form $\sw \in \Omega(C)$ is  \df{exact} if $\int_\rho \sw = 0$
for any loop $\rho$ in $C$. The vector space of exact 1-forms on $C$ is denoted by $ \Omega_0( C )$.
For a path from a leaf to another or a loop $\rho$ in $C$,  define the 1-form $\sw^\rho \in \Omega(C)$ \df{dual to} $\rho$ by
\[
\sw^\rho(\vec{e}) =
\left\lbrace
\begin{array}{rl}
1 \quad&  \text{if } \; \vec{e} \in \im(\rho) \\
-1 \quad&  \text{if }\; -\vec{e} \in \im(\rho) \\
0 \quad& \text{otherwise}
\end{array}
\right. .
\]
\end{definition}

\begin{proposition}\label{tdiff}
Let $C$ be a tropical curve of genus $g$ with $n\geqslant 2$ leaves.

-- a) The sum of the residues of any 1-form $\sw\in \Omega(C)$ is zero.

-- b) $\dim_\R \Omega_0 (C)= n-1$, $\dim_\R \Omega_\cH (C)= g$ and $\Omega (C) = \Omega_0 (C) \oplus \Omega_\cH (C)$.

-- c) Any element of $\Omega_0 (C) $ is determined by its residues.
\end{proposition}

\begin{proof} 
Let us choose $g$ edges $e_1,\dots, e_g$ in $C$ such that the graph $C^*$ obtained from cutting $C$ at the middle of the latter $g$ edges is a connected tree. For $1\leqslant j \leqslant g$, there is exactly one path (up to orientation) in $C^*$ joining the two ends given by the cut at $e_j$. This path induces a loop $\rho_j$ in $C$. Denote $\sw_j \in \Omega(C)$ the $1$-form dual to the latter loop. Given an appropriate orientation $\vec e_j$ on $e_j$, we have that $\int_{\vec e_j}\sw_k=\delta_{jk}$ where $\delta_{jk}$ is the Kronecker delta. Label $l_1,\dots,l_n$ the leaves of $C$. For $1\leqslant j \leqslant n-1$, there is only one path from $l_j$ to $l_n$ in $C$ coming from a path in $C^*$. Denote this path by $\rho_{g+j}$ and denote by $\sw_{g+j}$ its dual $1$-form, see Definition \ref{pathloop}. We claim that the elements $\sw_1,\dots,\sw_{g+n-1}$ in $\Omega(C)$ are linearly independent. Indeed, the matrix of the linear map from $\big( \bigoplus_j \sw_j \cdot\R  \big)$ to $\R^{g+n-1}$ given by $\big(\int_{\vec e_1} ,\dots,\int_{\vec e_g} , \Res_{ l_2} ,\dots,\Res_{ l_n}\big)$ is diagonal and invertible. Now, we claim that $\Omega(C)$ has dimension $g+n-1$ so that $\sw_1,\dots,\sw_{g+n-1}$ form a basis. Indeed, for any $\sw\in \Omega(C)$, there exists a unique linear combination $\widetilde \sw$ of $\sw_{g+1},\dots,\sw_{g+n-1}$ such that $\sw-\widetilde \sw$ has no residues and hence defines an element of $H^1(C,\R)$. From the previous arguments, the elements $\sw_1,\dots,\sw_{g}$ are linearly independent in the $g$-dimensional vector space $H^1(C,\R)$. It follows that $\Omega(C)\simeq H^1(C,\R)\oplus \big(\bigoplus_{j=g+1}^{g+n-1} \sw_j \cdot \R\big)$. In particular, the vector space $\Omega(C)$ has expected dimension. Now, since the sum of the residues of every element of the basis $\sw_1,\dots,\sw_{g+n-1}$ is zero, the point $a)$ holds. As shown previously, we have that $\Omega_\cH(C)\simeq H^1(C,\R)$ is $g$-dimensional. In particular, there exists elements $\left\lbrace \widetilde \sw_1,\dots,\widetilde\sw_g \right\rbrace$ in $\Omega_\cH(C)$ such that $\int_{\rho_k} \widetilde \sw_j= \delta_{jk}$ for any $1\leqslant k \leqslant g$. Thus, there exists for any $\sw\in \Omega(C)$ a unique linear combination $\widetilde \sw$ of $\widetilde \sw_1,\dots,\widetilde\sw_g$ such that $\sw-\widetilde \sw$ is exact. It follows that $\Omega (C) = \Omega_0 (C) \oplus \Omega_\cH (C)$ and $\dim_\R \Omega_0 (C)= n-1$, proving $b)$. For $c)$, consider the unique $\widetilde \sw_j$, $g<j\leqslant g+n-1$, such that $\sw_j-\widetilde \sw_j \in \Omega_0(C)$. By construction, the element $\widetilde \sw_j$ has residues $1$ at  $l_j$ and $-1$ at $l_{n}$ and no other residues. It follows that for any $\sw\in\Omega_0(C)$, the element $\sw -\big( (\Res_{l_1}\sw)\widetilde \sw_{g+1} +\dots + (\Res_{l_{n-1}}\sw) \widetilde \sw_{g+n-1}\big)$ has no residues. Since it is also exact by construction, it is zero by $b)$. This proves the point $c)$.
\end{proof}

\begin{definition}\label{defPir}
A \df{collection of residues} of dimension $m$ on a simple tropical curve $C$ with $n$ leaves is real $n \times m$ matrix $R:= \big( r_{k,j}  \big)_{k,j}$ such that $\sum_j r_{k,j} =0$ for any $1 \leqslant k \leqslant m$.
Denote by $$\sw_R^C := (\sw^{C}_{R,1},\dots, \sw^{C}_{R,m})$$ the vector of exact 1-forms on $C$ induced by the $m$ rows of $R$, see Proposition \ref{tdiff}. In practice, we will simply denote the latter vector by $\sw_R=(\sw_{R,1}, \dots, \sw_{R,m})$ when no confusion is possible. Finally, define the map
\[
\df{
\begin{array}{rcl}
\pi_R  \; : \; C & \rightarrow & \R^m\\
q & \mapsto & \big(  \int^q_{p} \sw_{R,1},\dots,  \int^q_{p} \sw_{R,m} \big)
\end{array}
}
\]
up to the choice of an initial point $p \in C$. 
\end{definition}

\begin{remark}
Exact 1-forms are exactly those forms obtained as gradient of particular functions on tropical curves. These functions are the harmonic tropical morphisms that we will introduce in Section \ref{secharmtrop}.
\end{remark}
\medskip

\subsubsection{Phase-tropical curves}

Phase-tropical curves first appeared in \cite{Mikh05} under the name of complex tropical curves. In the latter reference, Mikhalkin considered the diffeomorphism 
\[
\begin{array}{rcl}
H_t \; \;  : \;\; (\bC^\star)^2 & \rightarrow & (\bC^\star)^2 \\
(z,w) & \mapsto & \left( \vert z \vert^{\frac{1}{\log(t)}} \frac{z}{\vert z \vert} ,\vert w \vert^{\frac{1}{\log(t)}} \frac{w}{\vert w \vert} \right)
\end{array}
\]
satisfying $\A \circ  H_t = \frac{1}{\log(t)} \cdot \A$. Now,
for a family of algebraic curves $\left\lbrace \C_t \right\rbrace_{t>1} \subset (\mathbb{C}^\star)^2$ converging to a tropical curve $C\subset \R^2$ in the sense of \eqref{eq:conv}, the limit of the family  $\left\lbrace H_t(\C_t) \right\rbrace_{t>1} \subset (\mathbb{C}^\star)^2$ in Hausdorff distance on compact sets is a so-called phase-tropical curve $V\subset\ttor$. In particular, we have that $\A(V)=C$. 

For the purpose of this paper, we will introduce abstract phase-tropical curves. Such objects have already been considered in the unpublished work \cite{Mikhunpub} and we will see below that there are at least two different ways to think about them. The work \cite{Mikhunpub} gives yet an other point of view.

Consider the line $ \mathcal{L} := \left\lbrace (z,w) \in (\bC^\star)^2 \: \big| \: z+w+1=0  \right\rbrace$. The sequence of topological surfaces $\big\lbrace H_t(\mathcal{L}) \big\rbrace_{t >1}$ converges in Hausdorff distance to the phase-tropical line $L$ when $t$ tends to $\infty$. The subset $L\subset \ttor$ can be described as follows. The amoeba $\A(L)\subset \R^2$ is the tropical line $\Lambda$ consisting of $3$ half-rays emanating from $(0,0)$ and directed by $(-1,0)$, $(0,-1)$ and $(1,1)$ respectively. If $r$ is the open half-ray directed by $(-1,0)$ (respectively  $(0,-1)$, respectively $(1,1)$), the set $\A^{-1}(r) \cap L$ is given by the intersection of $\A^{-1}(r)$ with the cylinder $\{z=1\}$ (respectively $\{w=1\}$, respectively $\{z=-w\}$) and is therefore an open annulus. These $3$ annuli are glued to the preimage of $(0,0)\in \R^2$ under the map $\A_{\vert L}$. This preimage is the closure of the so-called \df{coamoeba} $\Arg(\mathcal{L}) \subset (S^1)^2$ where $\Arg$ is the component-wise argument map $\Arg : (\bC^\star)^2\rightarrow (S^1)^2$. The latter (closed) coamoeba is the topological pair of pants illustrated in Figure \ref{fig:coamoeba}: it consists in the union of 2 triangles delimited by 3 geodesics in the argument torus $\sone$. Therefore, the phase-tropical line $L$ is also a pair of pants, as illustrated in  Figure \ref{fig:L0}. We also refer to \cite[Proposition 6.11]{Mikh05}.

\begin{figure}[h]
\centering
\input{L02.ttex}
\caption{The map $ \A : L \rightarrow \Lambda$. The preimage of the vertex of $\Lambda$ in $L$ is pictured in grey.}
\label{fig:L0}
\end{figure}

The phase-topical line $L$ is globally invariant under complex conjugation, since each $H_t(\mathcal{L})$ is. The \df{fixed locus} $\R L\subset L$ consists of $3$ arcs passing respectively through one of the $3$ vertices of the coamoeba $\A^{-1}(0,0)\cap L$ and joining $2$ of the $3$ punctures of $L$. Observe that for any point $p\in \Lambda\setminus (0,0)$, the geodesic $\A_{\vert L}^{-1}(p) \subset \sone$ intersects $\R L$ in exactly two points.

In the present text, a \df{toric transformation} is a map  $ B : \ttor \rightarrow \ttor $ of the form 
\[ (z,w) \mapsto \big( z_0 z^{a}w^{b}, w_0 z^{c}w^{d} \big) \]
where $(z_0,w_0) \in  \ttor$ and $ \begin{psmallmatrix} a & b \\ c & d \end{psmallmatrix} \in \text{GL}_{2} \big( \mathbb{Z} \big)$.
The map $B$ descends to an affine linear transformation on $\mathbb{R}^2$ (respectively on $\sone$) by composition with the projection $\A$ (respectively $\Arg$) that we still denote by $B$. 
A \df{general phase-tropical line} in $\ttor$ is the image of $L$ by any toric transformation. For any generalised phase-tropical line $L'=B(L)$, we define $\R L':=B(\R L)$. Since the $6$ toric transformations mapping $L$ to itself fix $\R L$ globally, the latter definition does not depend on $B$.

We will now define abstract simple phase-tropical curves. A simple phase-tropical curve will be an abstract topological surface $V$ together with a map $\A_V : V \rightarrow C$ onto a simple tropical curve $C$. 
For the construction, we need the following data:

-- For any $v\in V(C)$, choose an injective continuous map $\pi_v: T_v \rightarrow \R^2$  satisfying the integrality and balancing conditions of Definition \ref{tropmorph} and choose also a generalised phase-tropical line $L_v \subset \ttor$ such that $\pi_v(T_v) \subset \A(L_v)$.

-- For any open, oriented edge $\vec e \in E(C)$ from a vertex $v$ to $v'$, choose a toric transformation $B_{\vec e} : \ttor \rightarrow \ttor $ such that $B_{\vec e}$ maps the annulus $\A^{-1}\big(\pi_v(e)\big)\cap L_v$ to the annulus $\A^{-1}\big(\pi_{v'}(e)\big)\cap L_{v'}$ and such that $(\pi_{v'})^{-1}\circ B_{\vec e}\circ (\pi_v)_{\vert e}=\id_e$. Define $B_{-\vec e}:= (B_{\vec e})^{-1}$.

Define $V$ to be the quotient of the space 
\[
 \coprod_{v \in V(C)}   \big\lbrace (p,q) \in L_v \times T_v \, \big| \, \A(p)= \pi_v(q) \big\rbrace 
\]
under the relation: $(p,q) \sim (p',q')$ $\Leftrightarrow$ $\big(q = q'$, $p\in L_v$, $p'\in L_{v'}$, $\exists \, \vec e \in E(C)$ from $v$ to $v'$ such that $B_{\vec e}(p)=p' \big)$.

Each \df{chart} $\big\lbrace (p,q) \in L_v \times T_v \, \big| \, \A(p)= \pi_v(q) \big\rbrace $ is in natural bijection with the topological pair of pants $\A^{-1}\big(\pi_v(T_v)\big)\cap L_v$ and equipped with the projection $\A_V:= (\pi_{v})^{-1}\circ \A$ onto $T_v$. It is an easy exercise to see that the above equivalence relation is compatible with the projections $\A_V$ defined on every chart. We then obtain a surjective map $\A_V: V\rightarrow C$ modelled on $\A:L_v \rightarrow \A(L_v)$ in a neighbourhood of each vertex $v\in V(C)$. In particular, the surface $V$ is a topological surface of genus $g=b_0(C)$ with $n$ punctures, where $n$ is the number of leaves of $C$.

\begin{figure}[h]
\centering
\scalebox{1}{
\input{coamoeba2.ttex}}
\caption{The coamoeba $ \Arg(\cL)$ (dark grey) and its 3 boundary geodesics (blue) in the fundamental domain $(-1,1]^2$ of $\sone$.}
\label{fig:coamoeba}
\end{figure}

\begin{definition}\label{def:phasetrop} 
A \df{simple phase-tropical curve} is a topological surface $V$ together with a map $\A_V:V\rightarrow C$ onto a simple tropical curve $C$ constructed from the data $\left\lbrace (\pi_v, L_v)\right\rbrace_{v\in V(C)}$, $\left\lbrace B_{\vec e} \right\rbrace_{e\in E(C)}$ as above.
\end{definition}

\begin{remark} For a family of algebraic curves $\left\lbrace \C_t \right\rbrace_{t>1} \subset (\mathbb{C}^\star)^2$ converging to a simple tropical curve $C\subset \R^2$ in the sense of \eqref{eq:conv}, the set $V:=\lim_{t\rightarrow \infty} H_t(\C_t) \subset (\mathbb{C}^\star)^2$ is a simple phase-tropical curve in the sense of Definition \ref{def:phasetrop}. The map $\A_V$ is given by $\A$, the maps $\pi_v$ are given by restriction of the embedding $C\hookrightarrow \R^2$ and the maps $B_{\vec e}$ are the identity.
\end{remark}

It is quite clear that any embedded phase-tropical curve $V\subset \ttor$ can be constructed from different data $\left\lbrace (\pi_v, L_v)\right\rbrace$, $\left\lbrace B_{\vec e} \right\rbrace$. Also, it seems natural to consider the phase-tropical curves $V\subset \ttor$ and $B(V)$ as isomorphic for any toric transformation $B:\ttor \rightarrow \ttor$.

\begin{definition}\label{def:isopt}
Two simple phase-tropical curves $A_V:V\rightarrow C$ and $A_{V'}:V'\rightarrow C'$ constructed respectively from the data $\left\lbrace (\pi_v, L_v)\right\rbrace$, $\left\lbrace B_{\vec e} \right\rbrace$ and $\left\lbrace (\pi'_v, L'_v)\right\rbrace$, $\left\lbrace B'_{\vec e} \right\rbrace$ are \df{isomorphic} if the following conditions are satisfied:

-- There exists an isomorphism $h:C\rightarrow C'$.

-- For any $v\in V(C)$, we denote by $B_v:\ttor \rightarrow \ttor$ the unique toric transformation such that $\pi'_{h(v)}=B_v \circ \pi_v$ and $B_v(L_v)=L'_{h(v)}$.

-- For any $\vec e\in E(C)$ from $v$ to $v'$, the restriction of $B'_{h(\vec e)}$ to $L'_{h(v)}$ is given by $B_{v'}\circ B_{\vec e} \circ (B_v)^{-1}$.
\end{definition}

Let us now show that simple phase-tropical curves can be described in a way that is similar to the description of Riemann surfaces using Fenchel-Nielsen coordinates.

Let $\A_V:V\rightarrow C$ be a simple phase-tropical curve. The underlying tropical curve $C$ can be described as a pair $(G, \ell)$ of a cubic graph $G$ equipped with a length function $\ell$, see Remark \ref{rem:tropGL}. The collection of maps $\pi_v$ of Definition \ref{def:phasetrop} induces a ribbon structure $\scrR$ on $G$ by pulling back the cyclical ordering of the leaves/edges of $\pi_v(T_v)$ given by the counter-clockwise orientation of $\R^2$. In order to define the twists along the edges of $G$, consider the following construction.

Choose an edge $e\in E(G)$ with adjacent vertices $v$ and $v'$. For any point $p \in e$, the chart $ \big\lbrace (p,q) \in L_v \times T_v \, \big| \, \A(p)= \pi_v(q) \big\rbrace$ allows to identify $\A_V^{-1}(p)$ to $S^1$ as follows. In the latter chart, the set $\A_V^{-1}(p)$ is a geodesic in the argument torus $\sone\subset \ttor$. In order to identify $\A_V^{-1}(p)$ with $S^1$, we need an orientation and the choice of an origin on the geodesic $\A_V^{-1}(p)\subset \sone$. We choose the orientation whose direct normal vector $\vec{n}$ is such that $T \A_V (\vec{n})$ is pointing towards $v$. For the origin, recall that there are exactly two points of $\R L_v$ on the geodesic $\A_V^{-1}(p)$. These two points belong to different components of $\R L_v$. The components of $\R L_v$ correspond to pair of leaves of $\A(L_v)$. Leaves are cyclically ordered by the counter-clockwise orientation of $\R^2$ and so are the pairs of them. In turn, the two points of $\A_V^{-1}(p)\cap \R L_v$ are ordered. We choose the first of these two points as the origin of $\A_V^{-1}(p)$. For instance, assume that $L_v=L$ and $p$ sits on the leaf of direction $(1,1)$ of $\Lambda$. Then $\A_V^{-1}(p)$ is the geodesic $\left\lbrace \arg(z)=-\arg(w)\right\rbrace$. The orientation chosen above is given by the vector $(-1,-1)$ and the origin is the point $(-1,1)\in \sone$.

The same construction can be repeated on the chart $ \big\lbrace (p,q) \in L_{v'} \times T_{v'}  \, \big| \, \A(p)= \pi_{v'} (q) \big\rbrace$ so that we get two identifications $\iota_v, \, \iota_{v'} : \A_V^{-1}(p)\rightarrow S^1$ satisfying
\[
\begin{array}{rcl}
\iota_{v'} \circ \iota_v^{-1} : S^1 & \rightarrow &  S^1\\
z & \mapsto & -\overline{\theta(e)z}.
\end{array}\]
for some $\theta(e)\in S^1$. The latter map is an involution. In particular, the element $\theta(e)$ does not depend on the ordering of $v$ and $v'$. Observe also that $\theta(e)$ does not depend on the choice of the point $p \in e$. We then obtain a map $\theta: E(G)\rightarrow S^1$.

\begin{proposition}\label{prop:ptFN}
The map $V\mapsto (G,\scrR, \ell, \theta)$ constructed above establishes a bijective correspondence between isomorphism classes of simple phase-tropical curves and quadruple $(G,\scrR, \ell, \theta)$ up to the equivalence relation generated by the following:

-- $(G,\scrR, \ell, \theta)\sim (G',\scrR', \ell', \theta')$ if there exists an isomorphism of graph $h:G\rightarrow G'$ such that $\scrR$, $\ell$ and $\theta$ are the respective pullbacks of $\scrR'$, $\ell'$ and $\theta'$ by $h$.

-- $(G,\scrR, \ell, \theta)\sim (G,\scrR', \ell, \theta')$ if $\scrR=\scrR'$ at all vertices except at $v$ and $\theta=\theta'$ at all edges except at the $3$ edges adjacent to $v$ where $\theta=-\theta'$.
\end{proposition}

\begin{proof}
Let us first show that the map $V\mapsto (G,\scrR, \ell, \theta)$ is surjective. Indeed, the pair $(G,\ell)$ is the simple tropical curve supporting $V$ and can therefore be chosen arbitrarily. For the coordinate $\scrR$, consider the following construction. For $V$ given by $\left\lbrace (\pi_v, L_v)\right\rbrace$, $\left\lbrace B_{\vec e} \right\rbrace$, a vertex $v_0\in V(C)$ and a toric transformation $B_0$, consider the phase-tropical curve $V'$ given by $\left\lbrace (\pi'_v, L'_v)\right\rbrace$, $\left\lbrace B'_{\vec e} \right\rbrace$ such that: for any vertex $v\neq v_0$, then $(\pi'_v, L'_v)=(\pi_v, L_v)$ and $(\pi'_{v_0}, L'_{v_0})=(B_0\circ\pi_v, B_0\circ L_v)$; for any edge $e$ not adjacent to $v$, then $B'_{\vec e}=B_{\vec e}$ and $B'_{\vec e}=B_{\vec e}\circ B_0^{-1}$ for any edge $\vec e$ adjacent to $v$ and oriented outwards $v$. The resulting phase-tropical curve $V'$ is isomorphic to $V$ by construction. If $B_0:\R^2\rightarrow\R^2$ is orientation preserving, then $V'$ maps to the same quadruple $(G,\scrR, \ell, \theta)$ as $V$. Otherwise, the curve $V'$ maps to $(G,\scrR', \ell, \theta')$ where $(\scrR',\theta')$ is related to $(\scrR,\theta)$ by in the second relation of the statement. This implies in particular that the map $V\mapsto (G,\scrR, \ell)$ is surjective. For the remaining coordinate $\theta$, consider the following construction. Start again from $V$ given by $\left\lbrace (\pi_v, L_v)\right\rbrace$, $\left\lbrace B_{\vec e} \right\rbrace$. Choose a vertex $v_0\in V(C)$, an edge $\vec e_0$ adjacent to $v$ and oriented towards $v$ and a toric transformation $B_0(z,w)=(e^{i a \nu}z,e^{i b \nu}w)$ such that $(a,b)\in \bZ^2$ is the slope of $\pi_{v_0}(\vec e)$. Consider now the phase-tropical curve $V'$ given by $\left\lbrace (\pi'_v, L'_v)\right\rbrace$, $\left\lbrace B'_{\vec e} \right\rbrace$ such that: for any vertex $v\neq v_0$, then $(\pi'_v, L'_v)=(\pi_v, L_v)$ and $(\pi'_{v_0}, L'_{v_0})=(B_0\circ\pi_v, B_0\circ L_v)$; for any edge $e$ not adjacent to $v$ and $e_0$, then $B'_{\vec e}=B_{\vec e}$ and $B'_{\vec e}=B_{\vec e}\circ B_0^{-1}$ for any edge $\vec e\neq \vec e_0$ adjacent to $v$ and oriented outwards $v$. Since $B_0$ fixes globally the annulus $\A^{-1}\big(\pi_{v_0}(e_0)\big)\cap L_{v_0}$, the data $\left\lbrace (\pi'_v, L'_v)\right\rbrace$, $\left\lbrace B'_{\vec e} \right\rbrace$ defines a phase-tropical curve $V'$ with twist function $\theta'$ such that $\theta'=\theta$ for all edges different from $e_0$ and $\theta'(e_0)=\theta(e_0)\cdot e^{i\nu}$. This proves the surjectivity of the map $V\mapsto (G,\scrR, \ell, \theta)$.

It remains to show that that $V$ is isomorphic to $V'$ if and only if the corresponding data $(G,\scrR, \ell, \theta)$ are equivalent to each other. To see this, observe that any isomorphism $V\rightarrow V'$ can be decomposed as the composition of an isomorphism on the underlying graph $G$ (taken care of by the first relation of the statement) and a composition of isomorphisms as constructed above for the surjectivity on the coordinate $\scrR$ (taken care of by the second relation of the statement). Hence, isomorphic phase-tropical curves are mapped to the same equivalence class. It remains to see that non-isomorphic curves $V$ and $V'$ are mapped to different equivalence classes of quadruple. If $V$ and $V'$ do not have isomorphic underlying tropical curves, then they are clearly mapped to different classes. If they do, we can assume with no  loss of generality that $V$ and $V'$ are constructed from the same data $\left\lbrace (\pi_v, L_v)\right\rbrace$ by applying an isomorphism on $V'$. In particular, $V$ and $V'$ map to the same triple $(G,\scrR,\ell)$. From there, it is an easy exercise to see that $V$ and $V'$ have necessarily different twist functions. 
\end{proof}

We denote by $V(G,\ell,\theta)$ the (isomorphism class of) simple phase-tropical curve given by the quadruple $(G,\scrR, \ell, \theta)$.

\section{Approximation of harmonic tropical curves}

\subsection{Convergence of imaginary normalised differentials}\label{sec:convind}

The purpose of this section is to study the limit of i.n.d. on families of algebraic curves that converge in the sense of Definition \ref{tropconv}. In general, it is complicated to determine the limits of linear systems on families of algebraic curves, as addressed for instance in \cite[\S 5]{GK} and references therein.

Recall that a vector of residues $R=(r_1,\dots,r_n)$ and a family $ \left\lbrace S_t \right\rbrace_{t>1} \subset \mgn$ give rise to the family $ \left\lbrace ( S_t , \omega_{R}) \right\rbrace_{t>1}\subset \thb_{g,n}$, see Definition \ref{defAr}. To any simple tropical curve $C$,  we associate the unique stable curve $S(C) \in \bmgn$ whose dual graph is the underlying cubic graph of $C$. Then, any $\sw \in \Omega (C)$ induces a generalised meromorphic differential on $S(C)$ as follows. Since each irreducible component of $S(C)$ is a sphere with $3$ marked points, it suffices to prescribe the residues of the restriction of the differential to each sphere. For a vertex $v\in V(C)$ and any adjacent $\vec{e}\in LE(C)$ oriented toward $v$, prescribe the residue at the marked point dual to $e$ on the component of $S(C)$ dual to $v$ to be $\sw(\vec{e})$. We still denote by $\sw$ the induced generalised meromorphic differential on $S(C)$ so that the pair $(S(C),\sw)$ is a point in $\thb_{g,n}$. Reciprocally, observe that any generalised meromorphic differential on $S(C)$ comes from a 1-form $\sw\in\Omega(C)$.
%
%
The main result of this section is the following.

\begin{mytheorem}\label{convdif}
Let $\left\lbrace S_t \right\rbrace_{t >1} \subset \mgn $ be a family converging to a simple tropical curve $C$. Then, for any collection of residues $R=(r_1,\dots,r_n)$, the family $\left\lbrace (S_t, \omega_{R}) \right\rbrace_{t >1} \subset \thb_{g,n}$  converges to the point $\big(S(C), \sw_{R}\big)$ defined above, where $\sw_R\in\Omega(C)$ is the 1-form of  Definition \ref{defPir}.
\end{mytheorem}

While proving Theorem \ref{convdif}, we will need some terminology. Recall from Definition \ref{tropconv} that $S_t := \fn_G(\ell_t, \theta_t)$ comes with a decomposition into  generalised pairs of pants $\left\lbrace Y_{v,t} \right\rbrace_{v\in V(G)}$ bounded by geodesics $\left\lbrace \gamma_{e,t} \right\rbrace_{e\in E(G)}$. 

\begin{definition}\label{assloop}
Let $\left\lbrace S_t \right\rbrace_{t >1} \subset \mgn$ be a family converging to a simple tropical curve $C$ and $\rho=\left\lbrace \vec{e}_1,\dots, \vec{e}_k \right\rbrace$ be any loop in $C$ (see Definition \ref{pathloop}). Define $\rho_t \subset S_t$ to be the piecewise-geodesic oriented loop satisfying:

-- $\rho_t \subset \cup_{1\leqslant j \leqslant k} Y_{v_j,t}$ where $v_j$ the initial vertex of the oriented edge $\vec{e}_j$.

-- $\rho_t  \cap \itr\big(Y_{v_j,t}\big)$ is the unique component of $\R Y_{v_j,t}$ joining $\gamma_{e_{j-1},t}$ to $\gamma_{e_{j},t}$ and oriented towards the latter (take $e_0:=e_k$).

-- $\rho_t  \cap \gamma_{e_j,t}$ is either a point or an oriented arc strictly contained in $\gamma_{v_j,t}$ whose direct normal vector points towards $Y_{v_j,t}$.
\end{definition}

\begin{remark}\label{rem:pathconv}
The family of loops $\rho_t \subset S_t$ in the above definition converges to a loop $\rho_\infty$ in the limiting stable $S(C)$. This follows from the fact that the fixed locus $\R Y_{v_j,t}$ converges to the real part of the corresponding irreducible component of $S(C)$. Recall that such a component is bi-holomorphic to $\bC P^1$ with marked points $-1$, $1$ and $\infty$ and with real part $\R P^1$.
\end{remark}

\begin{proof}[Proof of Theorem \ref{convdif}.]
Assume first that the sequence $\big\lbrace \int_{\gamma_{e,t}} \omega_{R}
\big\rbrace_{t>1}$ converges for any $e\in E(G)$ when $t$ tends to $\infty$. According to Section \ref{sec:thb}, the latter is equivalent to the convergence of the sequence $\left\lbrace (S_t,  \omega_{R}) \right\rbrace_{t>1}$ in $ \thb_{g,n}$. In particular, the limiting differential, denote it $\sw$, restricted to any irreducible component of $S(C)$ is a meromorphic differential whose set of poles is a subset of the $3$ marked points of the component. Let us show that the limit is $\big(S(C), \sw_{R}\big)$.

As a direct consequence of Definition \ref{tropconv}, the family $\left\lbrace S_t \right\rbrace_{t>1}$ converges to the stable curve $S(C)$. Now, fix a loop $\rho=\left\lbrace \vec{e}_1,\dots, \vec{e}_k \right\rbrace$ in $C$ and consider the corresponding loop $\rho_t \subset S_t$ of Definition \ref{assloop}. For simplicity, denote $\gamma_{j,t}:=\gamma_{e_j,t}$ and  $K_{j,t}\subset S_t$ the 
collar around $\gamma_{j,t}$. Orient $\gamma_{j,t}$ so that the direct normal vector points toward $Y_{v_j,t}$ where $v_j$ is the initial vertex of the oriented edge $\vec{e}_j$.  Then, there exists a unique bi-holomorphism $\psi_{j,t} : K_{j,t} \rightarrow A_{m(j,t)}$
such that:

-- $\psi_{j,t} \big( \rho_{t} \cap K_{j,t}  \big) $ is a transversal path in $ A_{m(j,t)}$,

-- $\psi_{j,t} \big( \rho_{t} \cap K_{j,t}  \big)  \cap \left\lbrace z \in A_{m(j,t)} \, \big| \: Re(z) =0 \right\rbrace = \left\lbrace 0 \right\rbrace$.\\
We claim that the family $(\psi_{j,t})_\ast \, \omega_{R}$ defined respectively on $ A_{m(j,t)}$ converges in the sense of Definition \ref{deflocconv}. Indeed, the family of collars $K_{j,t}$ converges to a neighbourhood $K_{j,\infty}$ of the node $q_j\in S(C)$ given by the vanishing cycle $\gamma_{j,t}$. The neighbourhood $K_{j,\infty}$ is actually the union of the two cusp at $q_j$ according to Lemma \ref{propcolcusp}. In particular, the family $\left\lbrace K_{j,t}\right\rbrace_{1<t\leqslant\infty}$ is just an other model for $\left\lbrace C_{\lambda} \right\rbrace_{0\leqslant \lambda<\lambda_0}$ for some $\lambda_0\leqslant 1$. Since $\left\lbrace (S_t,  \omega_{R}) \right\rbrace_{t>1}$ converges in $\thb_{g,n}$ by assumption, the restriction of $\omega_R$ to $K_{j,t}$ converges in the sense Definition \ref{deflocconv}. The claim follows.

We now claim that the integral of $\omega_R$ on $\rho_t \setminus \big( \bigcup_{1\leqslant j \leqslant k} \,  K_{j,t} \big)$ converges when $t$ tends to $\infty$ and is therefore a bounded function in $t$. Indeed, the family of loops $\rho_t$ converge to a loop $\rho_\infty \subset S(C)$ passing through the node of $S(C)$ dual to the edges $e_1,\dots,e_k$, see Remark \ref{rem:pathconv}. The latter nodes are the only possible poles of the limiting differential $\sw$ lying on $\rho_\infty$. The collars $K_{j,t}$ converge to a neighbourhood $K_{j,\infty}$ of the latter nodes, see Lemma \ref{propcolcusp}. It follows that 
\[ \lim_{t\rightarrow \infty} \; \textstyle \int_{\rho_t \setminus \big( \bigcup_{1\leqslant j \leqslant k} \,  K_{j,t} \big)} \omega_R \; = \;  \int_{\rho_\infty \setminus \big( \bigcup_{1\leqslant j \leqslant k} \,  K_{\infty,t} \big)} \sw  \; < \; \infty \] and the claim is proven. It follows in turn that 
\[ \int_{\rho_{t}} \omega_{R} \; =  \;  \sum_{j=1}^{k} \textstyle \big( \int_{\rho_{t} \cap K_{j,t}} \omega_{R} \big)\; + \; O(1).\]
Applying Proposition \ref{convdifloc} to the restriction of $\omega_R$ to each $K_{j,t}$ (with $\alpha=1$ and $l=\ell(e_j)$), we deduce from the above formula that 
\begin{equation}\label{lim}
\lim_{t \rightarrow \infty} \;  \frac{1}{\log(t)} \int_{\rho_t} \omega_{R} \; = \; \sum_{j=1}^{k} \ell(e_{j}) \Lambda_{j}
\end{equation}
where $\Lambda_{j} := \displaystyle \lim_{t \rightarrow \infty} \; \textstyle \frac{1}{2\pi i} \int_{\gamma_{j,t}} \omega_{R}$. Since $\omega_{R}$ is an i.n.d. on any $S_t$, it follows that $\int_{\gamma_{j,t}} \omega_{R}$ and $\int_{\rho_t} \omega_R$ are purely imaginary, for any $t$ and $j$. In particular, we have that $ \Lambda_{j} \in \mathbb{R} $ for any $j$. Considering the real part on both sides in \eqref{lim}, we obtain that 
\begin{equation}\label{lim2}
\sum_{j=1}^{k} \ell(e_{j}) \Lambda_{j} = 0. 
\end{equation}
Since the above formula holds for any loop $\rho=\left\lbrace \vec{e}_1,\dots, \vec{e}_k \right\rbrace$, it follows that the limiting differential $\sw$ comes from an exact 1-form on $C$ with residue vector $R$. According to Proposition \ref{tdiff}, the differential $\sw$ on $S(C)$ necessarily comes from $\sw_{R}\in \Omega(C)$. This proves the theorem in the case when $\big\lbrace \int_{\gamma_{e,t}} \omega_{R}
\big\rbrace_{t>1}$ converges for any $e\in E(C)$.

Assume now that the sequence $\big\lbrace \int_{\gamma_{e,t}} \omega_{R}
\big\rbrace_{t>1}$ is bounded uniformly in $t$ for any $e\in E(C)$. For any discrete infinite subset $I \subset (1,\infty)$ such that $\big\lbrace \int_{\gamma_{e,t}} \omega_{R}
\big\rbrace_{t\in I}$ converges for any $e\in E(C)$, it follows from the previous case that the the family $\left\lbrace (S_t, \omega_{R}) \right\rbrace_{t \in I} \subset \thb_{g,n}$  converges to the point $\big(S(C), \sw_{R}\big)$. As any converging subsequence converges to the same limit, the original sequence $\left\lbrace (S_t, \omega_{R}) \right\rbrace_{t >1} \subset \thb_{g,n}$  converges to $\big(S(C), \sw_{R}\big)$.

Assume finally that the sequence $\big\lbrace \int_{\gamma_{e,t}} \omega_{R}
\big\rbrace_{t>1}$ is unbounded for at least one edge $e\in E(C)$. Define $M_t:= \max \big\lbrace 1, \max \big\lbrace \vert \int_{\gamma_{e,t}} \omega_{R} \vert  \big\rbrace_{e\in E(C)} \big\rbrace$. Then,  there exists a discrete infinite subset $I \subset (1,\infty)$ such that  $\displaystyle \lim_{t\rightarrow \infty}  M_t  =\infty$. In particular, the sequences $\big\lbrace \int_{\gamma_{e,t}} \widetilde \omega_{R}
\big\rbrace_{t\in I}$, where $\widetilde \omega_{R}:= \omega_{R}/M_t$, are bounded for any $e\in E(C)$. Extracting a subsequence if necessary, we can therefore assume that $\big\lbrace \int_{\gamma_{e,t}} \widetilde \omega_{R}
\big\rbrace_{t\in I}$ converges for any $e\in E(C)$. Consequently, the family $\left\lbrace (S_t,  \widetilde \omega_{R}) \right\rbrace_{t\in I}$ converges to a point $\big(S(C), \widetilde \sw\big)\in \thb_{g,n}$. Then, we can apply the same reasoning as in the first case to show that $\widetilde \sw$, seen as a 1-form on $C$, is exact. By definition of $M_t$, the 1-form $\widetilde \sw \in \Omega_0(C)$ is non-zero and has no residues at the leaves of $C$ since $\widetilde \omega_{R}^{S_t}$ has residues $R/M_t$ and  $\lim M_t=\infty$. In other words, the 1-form $\widetilde \sw$ is both exact and holomorphic. It follows from Proposition \ref{tdiff}-$b)$ that $\widetilde \sw$ is zero. This is a contradiction. We deduce that all the sequences $\big\lbrace \int_{\gamma_{e,t}} \omega_{R}
\big\rbrace_{t>1}$ are necessarily bounded and the theorem follows from the two previous cases.
\end{proof}
\medskip

\subsection{Convergence of periods}\label{sec:cop}

The refined notion of phase-tropical convergence allows to obtain a finer result than Theorem \ref{convdif}.
\begin{mytheorem}\label{thm:convperiod}
Let $\left\lbrace S_t \right\rbrace_{t >1} \subset \mgn $ be a family converging to a simple phase-tropical curve $V(G,\ell,\theta)$ and $R=(r_1,\dots,r_n)$ be a collection of residues. 
%
%
For any loop $\rho=\left\lbrace\vec{e}_1,\dots,\vec{e}_k \right\rbrace$ in the simple tropical curve $C:=(G,\ell)$ and $\rho_t$ the associated loop in $S_t$, see Definition \ref{assloop}, we have 
\[  \lim_{t \rightarrow \infty} \; \int_{\rho_t} \omega_{R}^{S_t}  \; = \;  \sum_{j=1}^k \;  \log \big( \theta (e_j) \big) \cdot \sw_{R}^{C} \big( \vec{e}_j \big)\]
where the branch of $\log$ is chosen such that $\log  :  S^1 \rightarrow \left[ 0, 2\pi i \right) \subset \mathbb{C}$.
\end{mytheorem}

In order to prove Theorem \ref{thm:convperiod}, we will study imaginary normalised differentials on $S \in 
M_{0,3}$. Let $\omega$ be such a differential on $S$. We are interested in paths $\rho\subset S$ for which the restriction to $\rho$ of the real-valued differential $\Im(\omega)$ is identically zero. We denote by $h$ the hyperbolic metric on $S$ and  define $\Im(\omega)^\vee$ (respectively $\Re(\omega)^\vee$) to be the vector field dual to $\Im(\omega)$ (respectively $\Re(\omega)$) with respect to $h$. The vector field $\Im(\omega)^\vee$ is 
obtained from $\Re(\omega)^\vee$ by a rotation with angle $\pi/2$. As a consequence, a path $\rho\subset S$ is such that $ 
\Im(\omega)_{\vert \rho} \equiv 0$ if and only if $\rho$ is parallel to the vector field $\Re(\omega)^\vee$.

We identify $ S \simeq \bC P^1 \setminus \left\lbrace -1, 1, \infty \right\rbrace$ with afine coordinate $z$ and assume for a moment that none of the residues of $\omega$ is zero. If necessary, we replace $\omega$ with $-\omega$ and apply an automorphism of $S$ exchanging the punctures so that 
\begin{equation}\label{omega}
\omega = \left( \frac{\lambda_{+}}{z-1} + \frac{\lambda_{-}}{z+1}\right) dz 
\end{equation} 
with $\lambda_-, \lambda_+ >0$. We denote by $\widetilde S$ the real oriented blow-up of $S$ at $-1$, $1$ and $\infty$ and denote by $\gamma_-$, $\gamma_+$ and $\gamma_\infty$ the respective boundary components of $\widetilde S$. The vector field $\Re(\omega)^\vee$ does not extend to the boundary of $\widetilde S$ as its modulus gets arbitrarily large. However, its asymptotic direction is well-defined. Below, we denote by $\R S$ the fixed locus of the unique anti-holomorphic involution on $S$ fixing the puncture and denote by $\R \widetilde S$ its lift to $\widetilde S$.

\begin{lemma} For any i.n.d. $\omega$ on $S$ as in \eqref{omega}, we have the following:

-- The 3 connected components of $\mathbb{R} S$ are parallel to $\Re(\omega)^\vee$.

-- $\Re(\omega)^\vee$ is asymptotically orthogonal to $\gamma_{-}$, $\gamma_+$ and $\gamma_\infty$. It is oriented inward $\widetilde S$ at $\gamma_{-}$ and $\gamma_+$ and outward at $\gamma_\infty$.

-- For any point $p$ in $\gamma_{-}$ or $\gamma_+$ away from  $\mathbb{R} \widetilde S$, the flow line of $\Re(\omega)^\vee$ starting at $p$ ends in $\gamma_\infty$.

-- The segment $[-1,1]$ is the only arc in $\widetilde S$ joining $\gamma_{-}$ to $\gamma_+$ that is parallel to $\Re(\omega)^\vee$.
\end{lemma}

\begin{figure}[h]
\centering
\input{PhasePortrait.ttex}
\caption{Phase portrait of $\Re(\omega)^\vee$ on $\widetilde S$.}
\label{fig:phpt}
\end{figure}

\begin{proof}
For the first part, write $z=x+iy$. Then, for any $z\in \R S$, that is $z=x$, we have that $\Re(\omega)= \left( \frac{\lambda_{+}}{x-1} + \frac{\lambda_{-}}{x+1}\right) dx$ at $z$. The dual to $dx$ with respect to the metric $h$ is a multiple of $\partial x$ since $h$ is conformally equivalent to the Euclidean metric. It follows that $\mathbb{R} S$ i parallel to $\Re(\omega)^\vee$. 

For the second part,  we compute the asymptotic of $\Re(\omega)$ when $z$ tends to $1$ (the computation for $-1$ and $\infty$ are similar). We have that 
\[
\begin{array}{rcl}
\Re(\omega)  &\underset{z\rightarrow 1}{\sim} & \displaystyle \Re\left( \frac{\lambda_{+}}{z-1} dz\right) =  \Re\left( \frac{\lambda_{+}}{z-1} \right)dx - \Im \left( \frac{\lambda_{+}}{z-1} dz\right) dy  \\ & & \\
& = & \displaystyle  \frac{\lambda_+}{\vert z-1\vert^2}\big( \Re(\bar{z}-1) dx - \Im(\bar{z}-1) dy\big) = \frac{\lambda_+}{t\big(\tilde x^2+ \tilde y ^2\big)}\big( \tilde x \, dx+\tilde y\,dy\big)
\end{array}
\]
where we wrote $z=1+t(\tilde x +i \tilde y)$. It follows that if $z$ tends to $1$ with asymptotic direction $(\tilde x, \tilde y)$, the differential $\Re(\omega)^\vee$ is asymptotically a positive multiple of $\tilde x \, \partial x+\tilde y \,\partial y$. The second part is proven.

For the last two parts, observe that $\Re(\omega)^\vee$ has a unique zero at the point $ \zeta := \frac{\lambda_+ - \lambda_{-}}{\lambda_+ + \lambda_{-}} \in (-1,1)$. Indeed, the point $\zeta$ is the unique zero of $\omega$ and hence a zero of $\Re(\omega)$. Reciprocally, a zero of $\Re(\omega)$ is also a zero of $\Im(\omega)$ since the two vector fields are obtained from each other by a rotation. Hence, a zero of $\Re(\omega)$ is also a zero of $\omega$. Now, the vector field $\Re(\omega)^\vee$ is the gradient vector field of the harmonic function $z\mapsto \int^z \Re(\omega)$ with critical point at $\zeta$. By the maximum principle and the Morse Lemma, the gradient field of the latter function is the gradient field of the function $(x,y)\mapsto x^2-y^2$ in suitable local coordinates centred at $\zeta$. It follows that there are only $2$ incoming flow lines of $\Re(\omega)^\vee$ at $\zeta$ and $2$ outgoing flow lines, any other flow line does not meet $\zeta$. We deduce from the previous points that the incoming flow lines are the segment $(-1, \zeta)$ and $(\zeta,1)$. It implies that any flow line starting at a point in $\big(\gamma_+ \cup \gamma_-\big) \setminus \R \widetilde S$ avoids the only singular point $\zeta$ of  $\Re(\omega)^\vee$ and ends therefore at the boundary of $\widetilde S$. The remaining parts of the lemma follow now from the first two parts.
\end{proof}

\begin{lemma}
Let $\left\lbrace (Y_t,\omega_t)\right\rbrace_{t >1}$ be a family of pair of pants $Y_t$ equipped with a meromorphic differential $\omega_t$ whose only poles are at the possible punctures of $Y_t$. Assume that $Y_t$ converges to $S$ and that $\omega_t$ converges to an i.n.d. $\sw$ on $S$. For any connected component $\rho_t \subset Y_t$ of $\mathbb{R} Y_t$, we have 
\[ \lim_{t \rightarrow \infty} \int_{\rho_t} \Im(\omega_t) = 0. \]
\end{lemma}

\begin{proof}
Assume first that the differential $\sw$ is as in \eqref{omega}. For any $t>1$, there exist $r_{t,-}, \, r_{t,+}, \, r_{t, \infty} >0$ such that 
\begin{equation}\label{eq:model}
Y_t \simeq \left\lbrace z \in \mathbb{C} \; \big| \; \vert z \vert \leqslant r_{t,\infty}, \; \vert z-1 \vert \geq r_{t,+}, \; \vert z+1 \vert \geq r_{t,-} \right\rbrace
\end{equation}
with boundary components $\gamma_{t,-}$, $\gamma_{t,+}$, and $\gamma_{t,\infty}$. Let us prove the lemma for $\rho_t = [r_{t,-}, r_{t,+}]$ (the two remaining cases are similar). According to the previous lemma, there exists $N>0$ such that for almost all $t>N$, we have  

-- $\int_{\gamma_{t,-}} \Im(\omega_t) >0$ and $\int_{\gamma_{t,+}} \Im(\omega_t)  >0$,

-- $\Re(\omega_t)^\vee$ is everywhere transversal to $\gamma_{t,-}$,  $\gamma_{t,+}$ and $\gamma_{t,\infty}$,

-- the differential $\omega_t$ has a single zero on $Y_t$.\\
Applying the same reasoning as in the proof of the previous lemma, we deduce the existence of $M>N$ such that for almost all $t>M$, there exists a unique arc $\breve{\rho}_t \subset Y_t$ parallel to $\Re(\omega_t)$ and joining $\gamma_{t,-}$ to $\gamma_{t,+}$. Necessarily, the arc $\breve{\rho}_t$ is in an $\varepsilon_t$-neighbourhood of $\rho_t$ with $\displaystyle \lim_{t\rightarrow \infty} \varepsilon_t=0$. The segment $\rho_t$ is homotopic to $\bar{\rho}_t:= \rho_{t,+} \circ \breve{\rho}_t \circ \rho_{t,-} $ where $ \rho_{t,-} $ is an arc in $\gamma_{t,-}$ and $ \rho_{t,+} $ is an arc in $\gamma_{t,+}$. As $\breve{\rho}_t$ converges to $\rho_t$, the arcs $\rho_{t,-}$ and $\rho_{t,+}$ eventually shrink to a point. It follows that
\[\lim_{t \rightarrow \infty} \int_{\rho_t} \Im(\omega_t) \;=\; \lim_{t \rightarrow \infty} \int_{\bar \rho_t} \Im(\omega_t) \;=\; \lim_{t \rightarrow \infty} \int_{\breve{\rho}_t} \Im(\omega_t) \;=\; 0 \]
and the lemma is proven in this case.

For a general sequence $\left\lbrace \omega_t \right\rbrace_{t >1}$, we can consider the sequence $ \omega_t'':=\omega_t+\omega_t'$ where $\omega_t'$ is the restriction of the form 
\[ \left( \frac{M}{z-1} + \frac{M}{z+1}\right) dz\]
to $Y_t$. For $M>0$ arbitrarily large, the sequences $\left\lbrace \omega_t' \right\rbrace_{t >1}$ and $\left\lbrace \omega_t'' \right\rbrace_{t >1}$ are as in the above paragraph. Hence, they satisfy the conclusion of the lemma and so does $\left\lbrace \omega_t \right\rbrace_{t >1}$.
\end{proof}

In the course of the proof of Theorem \ref{thm:convperiod} and later on in the text, we will need the following elementary observation.

\begin{remark}\label{rem:blup}
Given a family $Y_t$ as in \eqref{eq:model}, it is useful to compare two different types of metric on the boundary components of $Y_t$. On the one hand, we can consider on a given boundary geodesic $\gamma\subset Y_t$ the metric of length $2\pi$ obtained by rescaling the hyperbolic metric with $\frac{2\pi}{l(\gamma)}$. On the other hand, we can consider on the  boundary geodesic $\lambda_{t,+}$ (respectively $\lambda_{t,-}$, respectively $\lambda_{t,\infty}$) the pullback of the Euclidean metric on $S^1$ by the function $\arg(z-1)$ (respectively $\arg(z+1)$, respectively $\arg(1/z)$). If for given $t$ and $\gamma\subset Y_t$, the two above metrics might differ, the induced metrics on the limiting boundary component of $\widetilde S$ will coincide at $t=\infty$. In particular, twist parameters of families $\left\lbrace S_t \right\rbrace_{t>1}$ converging tropically may be interpreted as angles when $t$ tends to $\infty$.
\end{remark}

\begin{proof}[Proof of Theorem \ref{thm:convperiod}.]
For any loop $\rho \subset C$, the associated loop $\rho_t \subset S_t$ is made out of connected components of $\mathbb{R} Y_{v,t}$ for all the vertices $v\in V(C)$ in $\rho$ and arcs in the geodesics $\gamma_{e,t}$ for all the edges $\vec e$ in $\rho$, see Definition \ref{assloop}. By the previous lemma, the  parts of $\rho_t$ contained in the different $\mathbb{R} Y_{v,t}$ do not contribute to the limit of the integral $\int_{\rho_t}\Im(\omega_t)$. In order to prove the theorem, it suffices then to show that for any $\vec e \in \rho$, we have
\[  \lim_{t \rightarrow \infty} \int_{\rho_t \cap \gamma_{e,t}}  \Im(\omega_R) \; = \; \log \big( \theta(e) \big) \cdot \sw_{R} (\vec e). \]
For $t$ large enough, the flow lines of $\Re(\omega_t)^\vee$ are transversal to $\gamma_{e,t}$ and points toward the pair of pants $Y_{v,t}$ for a vertex $v$ adjacent to $e$. Identify the latter pair of pants to a subset of $\bC \setminus \{-1,1\}$ as in \eqref{eq:model} such that $\gamma_{e,t}$ maps to  $\gamma_{t,+}$. Then, for $\varepsilon >0$ small enough, the flow lines of $\Re(\omega_t)^\vee$ are also transversal to the boundary of the $\varepsilon$-neighbourhood $U_{\varepsilon,t} \subset Y_{v,t}$ with respect to the Euclidean distance in $\bC \setminus \{-1,1\}$. Denote $c_{t,1}$ and $c_{t,2}$ the flow lines of $\Re(\omega_t)^\vee$ starting at the end points of $\rho_t \cap \gamma_{e,t}$ and ending on $\partial U_{\varepsilon,t}$. Denote by $\rho_{t,\varepsilon} $ the arc on $\partial U_{\varepsilon,t} $ between the end points of $c_{t,1}$ and $c_{t,2}$ such that $\rho_t \cap \gamma_{e,t}$ is homotopic to $ \bar\rho_{t, \varepsilon}:=c_{t,2}^{-1} \circ \rho_{t, \varepsilon} \circ c_{t,1}$ (or $c_{t,2} \circ \rho_{t, \varepsilon} \circ c_{t,1}^{-1}$ depending on how $c_{t,1}$ and $c_{t,2}$ are oriented). For any arbitrarily large $t$ and small $\varepsilon$, we have
\[  \int_{\rho_t \cap \gamma_{e,t}}  \Im(\omega_R) = \int_{\bar\rho_{t, \varepsilon}}  \Im(\omega_R) = \int_{\rho_{t, \varepsilon}}  \Im(\omega_R) \]
since $c_{t,1}$ and $c_{t,2}$ are flow lines of $\Re(\omega_t)^\vee$. By Theorem \ref{convdif}, the differential $\omega_R$ converges to the differential $\sw_R$ on $Y_{v,\infty}$ when $t$ tends to $\infty$. Since the twist parameter $\theta_t(e)$ converges to $\theta(e)$ by assumption, the arc $\rho_{\infty, \varepsilon}\subset Y_{v,\infty}$ is joined by flow lines of $\Re(\sw_R)^\vee$ to an arc $\alpha:=\left\lbrace z \in S^1 \, \big| \, \theta \leqslant z \leqslant \theta \cdot \theta(e) \right\rbrace$ inside the boundary component $\gamma_+ \simeq S^1$ of $\widetilde S$, see Remark \ref{rem:blup}. It follows that 
\[  
\begin{array}{rcl}
\displaystyle \lim_{t \rightarrow \infty} \int_{\rho_t \cap \gamma_{v,t}}  \Im(\omega_R) &= &
\displaystyle \lim_{\varepsilon \rightarrow 0} \lim_{t \rightarrow \infty} \int_{\rho_{t, \varepsilon}}  \Im(\omega_{R} )   \\
& & \\
&=& \displaystyle \int_{z\in \alpha} \frac{\sw_{R} \big(\vec e)}{z}dz \; = \;\log \big( \theta(e) \big) \cdot \sw_{R}(\vec e).
\end{array}\]
\end{proof}

\begin{corollary}\label{cor:period}
Let $\left\lbrace S_t \right\rbrace_{t >1} \subset \mgn $ be a family converging to a simple phase-tropical curve $V:=V(G,\ell,\theta)$ and $R=(r_1,\dots,r_n)$ be a collection of residues. Then, the family of period matrices $\cP_{R,S_t}$ converges to a matrix $\cP_{R,V}$ depending only on $\theta$ and the projective class of $\ell$, that is $\cP_{R,(G,\ell,\theta)}=\cP_{R,(G,\lambda\cdot\ell,\theta)}$ for any $\lambda>0$.
\end{corollary}

\begin{proof}
Denote $C:=(G,\ell)$. According to Theorem \ref{convdif}, the period $ \frac{1}{2 i \pi} \int_{\gamma_{e,t}} \omega_{R}^{S_t}$ tends to $\sw_{R}^{C} (\vec e)$, where $\gamma_{e,t}$ is the geodesic of the pair of pants decomposition of $S_t$ corresponding to $e \in E(C)$. 
For any basis $\rho_1,\dots, \rho_g \in H_1(C, \mathbb{Z})$ and associated loops $\rho_{1,t},\dots, \rho_{g,t} \subset S_t$, Theorem \ref{thm:convperiod} implies that the period $\frac{1}{2 i \pi} \int_{\rho_{k,t}} \omega_{R,j}^{S_t} $ converges to a linear combination of the $\sw_{R}^{C} (\vec e)$, $e\in E(G)$, whose coefficients are determined by $\theta$. In particular, all the periods computed above depend only on the projective class of $\ell$ since $\sw_R^C$ does. As any period of $\omega_{R}^{S_t}$ can be computed from the above periods, the result follows.
\end{proof}

\begin{remark}\label{rem:jacobian}
The techniques used to prove Theorems \ref{convdif} and \ref{thm:convperiod} can be adapted to the compact case $n=0$ and $g\geqslant 2$ to compute the asymptotic of the period matrix of holomorphic differentials on families of algebraic curves. For a family $\left\lbrace S_t \right\rbrace_{t>1} \subset M_g$ converging to a simple phase-tropical curve $V(G,\ell,\theta)$, consider the geometric symplectic basis $\alpha_{1,t},\dots,\alpha_{g,t}, \beta_{1,t},\dots, \beta_{g,t}$ of $S_t$ (see \cite[\S 6.1.2]{farbmarg}) constructed as follows: the $\alpha_{j,t}$ consist of cycles of the form $\gamma_{e_j,t}$ where the $e_j\in E(G)$ are such that $G\setminus \big( \cup_j \, e_j \big)$ is a connected tree; the cycle $\beta_{j,t}$ is given by $\rho_{j,t}$ (see Definition \ref{assloop}) where $\rho_j$ is the unique cycle in $G\setminus \big(\cup_{k\neq j}\, e_k\big)$.  For the basis $\omega_{1,t},\dots,\omega_{g,t}$ of $\Omega^1(S_t)$ dual to $\alpha_{1,t},\dots,\alpha_{g,t}$, that is such that $\int_{\alpha_{j,t}} \omega_{k,t}=\delta_{jk}$, the matrix of periods $B_t$ defined by $(B_t)_{jk}:=\int_{\beta_{j,t}} \omega_{k,t}$ determines the Jacobian $J(S_t)$. We claim that 
\[ B_t \;\;  \underset{t\rightarrow \infty}{\sim} \; \; \log(t)\cdot B_{\Re} \, +\, i \cdot B_{\Im}\]
where $B_{\Re}$ is the matrix of the quadratic form $Q$ on $H_1(G,\bZ)$ defined in \cite[\S 6.1]{MZ} and $B_{\Im}$ is the matrix of the form $Q'$ such that $Q'(\rho_j)=-i \cdot \sum_{\vec e \in \rho_j} \log\big(\theta(e)\big)\in \R$.
%
%
We deduce in particular that the rescaled period matrix $\frac{1}{\log(t)}\cdot B_t$ converges to the tropical period matrix $Q$ of $C=(G,\ell)$. Moreover, the knowledge of the asymptotic of the imaginary part of $B_t$ suggests a phase-tropical compactification of the moduli space $A_g$ of principally polarised Abelian varieties, see \cite{BMV} and \cite{O2} for related matters. 
\end{remark}
\medskip

\subsection{Harmonic tropical morphisms}\label{secharmtrop}

For any planar curve $ \mathcal{C} \subset \ttor$, the authors of \cite{PR} introduced the spine of the amoeba $\A(\C)$, a canonical tropical curve contained in $\A(\C)$ such that $\A(\C)$ deformation-retracts on it.  A similar construction can be carried to the case of harmonic amoebas in $\R^2$, see \cite{Kri}. In this more general case, the spine is still a piecewise linear graph whose edges have now arbitrarily slopes, with the same deformation-retraction property. Such spines are part of a wider class of tropical curves that we introduce below.

\begin{definition}\label{htropmorph}
A \df{ harmonic tropical morphism} $ \pi \, : \, C \rightarrow \mathbb{R}^m $ on a tropical curve $C$ is a proper continuous map subject to the following 

--\textbf{Linearity}: for any $ e \in LE( C )$, the map $ \pi_{\vert e} $ is affine linear with respect to the metric on $e$.

--\textbf{Balancing}: For any $ v \in V ( C )$, denote by $\vec{e}_1, \dots, \vec{e}_n$ the  outgoing unitary tangent vectors to the $n$ leaves/edges adjacent to $v$. Then
\[ \sum_{1\leqslant j \leqslant n} d\pi (\vec{e}_{j}) \; = \; 0. \]
A \df{harmonic tropical curve} in $\mathbb{R}^m$ is the image of a harmonic morphism.
\end{definition}

Similarly to the case of harmonic amoebas, every harmonic tropical morphism can be described by integration of 1-forms, as stated below.

\begin{proposition}\label{propharmmorph}
Let $C$ be a simple tropical curve and $R$ a collection of residues of dimension $m$ on $C$. Then, the map $ \pi_R  :  C  \rightarrow  \mathbb{R}^m$ of Definition \ref{defPir} is a harmonic tropical morphism. For any vertex $v \in V(C)$ and any outgoing adjacent leaf/edge $\vec e$, the corresponding unitary tangent vector is mapped to $\omega_{R}^{C}(\vec e)$. 
Reciprocally, for any harmonic morphism $ \pi  \; : \; C  \rightarrow  \mathbb{R}^m$, there exists a unique collection of residues $R$ of dimension $m$ on $C$ such that $ \pi = \pi_R$.
\end{proposition}

\begin{proof}
Since the 1-form $\sw_{R}^{C}$ is a constant on any leaf/edge $e\in LE(C)$, the map $\pi_R$ is affine linear on any $e\in LE(C)$. By definition, the vector $\sw_{R}^{C} (\vec e)$ is the gradient of $\pi_R$ along the oriented edge $\vec e$.  The balancing condition of \ref{tropmorph} follows from the definition \ref{def1form}. The first part of the statement is proven.

Reciprocally, it is now clear that the map $\pi$ is given by integration of an $m$-tuple of constant 1-forms on any $e \in LE(C)$. The balancing condition of \ref{tropmorph} is clearly equivalent to the balancing condition of \ref{def1form}. Exactness follows from the fact that this $m$-tuple is given by differentiating the $m$-tuple of  coordinate functions of $\pi$. Uniqueness follows from Proposition \ref{tdiff}-$c)$.
\end{proof}
\begin{remark}
As for harmonic amoeba maps, the space of deformation of harmonic morphisms on a fixed tropical curve $C$ in $\mathbb{R}^m$ corresponds to the space of collection of residues of dimension $m$ on $C$, see Proposition \ref{tdiff}. It is then a real vector space of dimension $m(n-1)$, where $n $ is the number of leaves of $C$.
\end{remark}
\medskip

\subsection{Degeneration of harmonic amoebas}\label{sec:dha}

In order to prove theorem \ref{approxtrop}, we need to make precise the notion of convergence we use. Recall the notations introduced in Section \ref{sectropcurvmorph}.

\begin{definition}\label{GHconv}
The sequence of maps $ \frac{1}{\log (t)} \mathcal{A}_R : S_t \rightarrow \mathbb{R}^m $ \df{converges} to the tropical morphism $\pi_R : C \rightarrow \mathbb{R}^m$ if for any $v \in V(C)$ (respectively any $e\in E(C)$), the image of $Y_{v,t}$ (respectively $\gamma_{e,t}$) under the map $ \frac{1}{\log (t)}  \mathcal{A}_R (Y_{v,t}) $ converges in Hausdorff distance to $\pi_R(\htr_v)$ (respectively to the middle point of $\pi_R(e)$).
\end{definition}

Note that the latter notion of convergence is stronger than the convergence of the image of $ \frac{1}{\log (t)} \mathcal{A}_R $ to the image of $\pi_R$ in Hausdorff distance.

\begin{remark}\label{rem:odaka}
In \cite{O}, Odaka studied the convergence of families of compact Riemann surfaces with respect to the Gromov-Hausdorff distance. When each surface in the family is equipped with its hyperbolic metric  rescaled so that the diameter of the surface is $1$, Odaka showed that the limit of such a family is either a Riemann surface or a tropical curve, see \cite[Theorem  2.4]{O}. Although the idea of exhibiting tropical curves as limits of families of Riemann surfaces is common to \cite{O} and the present work, the two approaches are different.  Indeed, for a family of compact Riemann surfaces $\left\lbrace S_t\right\rbrace_{t>1} \subset M_g$ converging to $C:=(G,\ell)$ in the sense of Definition \ref{tropconv}, the width of the collar of any vanishing cycle of the family is asymptotically equivalent to $\log(\log(t))$. Therefore, the diameter of $S_t$ is equivalent to $D\cdot \log(\log(t))$ where $D$ is the combinatorial diameter of the graph $G$. In particular, the family $\left\lbrace S_t\right\rbrace_{t>1}$ converges in the sense of \cite{O} to the graph $G$ equipped with length $1/D$ on every edge.
\end{remark}

\begin{proof}[Proof of Theorem \ref{approxtrop}.]
Recall that the maps $\pi_R$ and $\A_R$ are defined up to the choice of an initial point for integration, see Definitions \ref{defAr} and \ref{defPir}. Let us choose the initial point of $\pi_R$ to be a vertex $v$ and  for $\A_R$, take any sequence of initial points $p_t\in S_t$ such that $p_t\in Y_{v,t}$ and $p_t$ sits in the complement of the half-collars and cusps of $Y_{v,t}$. For this choice of initial points, we will show below that the image of $Y_{v,t}$ (respectively $\gamma_{e,t}$) under $ \frac{1}{\log (t)}  \mathcal{A}_R (Y_{v,t}) $ converges in Hausdorff distance to $\pi_R(\htr_v)$ (respectively to the middle point of $\pi_R(e)$, for any $e\in E(C)$ adjacent to $v$). Since the integrals defining $\pi_R$ and $\A_R$ are additive with respect to the subdivision of the paths of integration, the above claim implies the statement.

Denote by $R_{v,t}$ the $m \times 3$ matrix of rescaled periods $\frac{1}{2\pi i} \int \omega_{R}^{S_t}$ along the 3 boundary geodesics of $Y_{v,t}$, oriented so that the direct normal vector is pointing inwards. The matrix $R_{v,t}$ is a matrix of residues of dimension $m$ that converges, thanks to Theorem \ref{convdif}, to a matrix $R_v$, columns of which correspond to the slope of $\pi_R$ on the leaves/edges of $\htr_v$. Recall that $S:=\bC P^1 \setminus \{-1,1,\infty\}$ and $Y_{v,t}$ can be identified with a subset of $S$ as in \eqref{eq:model}. In particular, the pair of pants $Y_{v,t}$ converges to the whole $S$ when $t$ tends to $\infty$ and the restriction of $\omega_{R}^{S_t}$ to $Y_{v,t}$ converges to $\omega_{R_v}^{S}$, see again Theorem \ref{convdif}. 
It follows that $ \mathcal{A}_R (Y_{v,t}) $ converges in Hausdorff distance to  $\mathcal{A}_{R_{v}} (S)$ when $t$ tends to $\infty$ (observe that the latter amoeba might be degenerate to a line or a point if $R_v$ has rank $1$, respectively $0$). 
As a direct consequence, we have that 
\[  \Gamma_v \; := \; \lim_{t \rightarrow \infty} \;  \frac{1}{\log(t)} \mathcal{A}_R (Y_{v,t}) \; \subset  \; \lim_{t \rightarrow \infty}  \; \frac{1}{\log(t)} \mathcal{A}_{R_{v}} (S) \; =: \; \Gamma_v^\infty .\]
The right-hand side $\Gamma^\infty_v$ is the only tripod with infinite edges (respectively the only line, respectively the only point) containing $\pi_R(\htr_v)$ if $R_v$ has rank $2$ (respectively $1$, respectively $0$). We claim that the left-hand side $\Gamma_v$ is $\pi_R(\htr_v)$. First, observe that if $\htr_v$ contains a leaf $e$ of $C$, then $\Gamma_v$ contains $\pi_R (e)$. Indeed, the set $\pi_R (e)$ is an infinite half-ray whose slope is given by the corresponding column of $R$ while the harmonic amoeba $\mathcal{A}_R (Y_{v,t})$ has a tentacle with the same slope. Consider now a bounded edge $e$ of $\htr_v$ such that the corresponding column in $R_v$ is non-zero, otherwise there is nothing to prove. Pick a point $s_t$ in the boundary component $\gamma_{e,t}$ of $Y_{v,t}$ for any $t>1$.
Then, the image of $s_t$ under $\frac{1}{\log(t)} \mathcal{A}_R$ is given by integrating $\frac{1}{\log(t)} \omega_{R}^{S_t}$ along any path $\rho_t$ from $p_t$ to $s_t$. In particular, the path $\rho_t$ has to go across the half-collar $\hc_{\gamma_{e,t}}$. Obviously, the sequence of paths $\{ \rho_t\}_{t>1}$ can be chosen so that their restriction to the half-collars form an admissible sequence. Therefore, Proposition \ref{convdifloc} applies with $\alpha=1/2$ and $l=\ell(e)$, since the modulus of the collar $\hc_{\gamma_{e,t}}$ is half the modulus of $K_{\gamma_{e,t}}$ which is $m\big(\ell_t(e)\big)$. It follows that the limit of $\frac{1}{\log(t)} \mathcal{A}_R (s_t)$ is the end point of the edge $\pi_R(e)$. Such points has to be the end points of $\Gamma_v$ by the maximum principle. This proves the theorem.
\end{proof}
\begin{remark}
Approximation of tropical curves by families of amoebas is treated for instance in \cite{Mikh05}, \cite{Mikh06}, \cite{NS}, or \cite{Nish}.
Considering harmonic amoebas allows to extend the various approximation theorems given in the above references, as shown in Theorem \ref{approxtrop}. Not only we can approximate harmonic tropical curves whose edges may have irrational slopes, but we can also approximate honest tropical curves that can not be obtained as limit of families of amoebas of algebraic curves, see Examples \ref{ex:cubic} and \ref{ex:cubicg2}.
\end{remark}

\begin{example}\label{ex:cubic}
Consider the tropical curve $C\subset \R^3$ pictured in Figure \ref{fig:tropcubic}. The latter is known not to be approximable by families of amoebas of algebraic curves, see \cite[Example 5.12]{Mikh04}. 
As a metric graph, all the edges of $C$ have length $1$ except the edge $e$ which has length $2$. Now, consider a family of curves $\{S_t\}_{t>1} \subset M_{1,12}$ decomposed into pairs of pants as shown on the right-hand side of Figure \ref{fig:tropcubic}. Require that every geodesic of the latter decomposition has length $\ell_t=\frac{2\pi^2}{\log(t)}$ except for the geodesic $\gamma$ which should have length $\frac{\ell_t}{2}$. Consider then the collection of residues $R$ of dimension $3$ on $S_t$ such that the gray, black, yellow and green punctures are respectively assigned the residue vectors $(0,-1,0)$,  $(-1,0,0)$, $(0,0,-1)$ and $(1,1,1)$. Theorem \ref{approxtrop} implies that $\frac{1}{\log(t)}\A_R(S_t)$ converges in Hausdorff distance to $C$. 
\end{example}

\begin{figure}[h]
\centering
\scalebox{0.9}{\input{tropcubic.ttex}}
\caption{Approximation of the tropical curve $C\subset \R^3$ by harmonic amoebas.}
\label{fig:tropcubic}
\end{figure}

\begin{example}\label{ex:cubicg2}
It is a classical fact that the geometric genus of an algebraic curve of degree $d$ in $\bC P^n$ is at most $\frac{(d-1)(d-2)}{2}$. This turns out not to be the case in tropical geometry, as shown in \cite{BBM2}.
Consider for instance the tropical curve $C\subset \R^3$ of degree $3$ and genus $2$ given in Figure 2 of the latter reference. As in the previous example, it is a simple exercise to construct a family $\{S_t\}_{t>1} \subset M_{2,12}$ and a collection of residue $R$ of dimension $3$ such that $\frac{1}{\log(t)}\A_R(S_t)$ converges in Hausdorff distance to $C$. 
\end{example}

\section{Approximation of phase-tropical curves}\label{sec:approxphasetrop}

In this last section, we come back to the algebraic case. As illustrated by the above examples, algebraicity impose severe restrictions on the realisability (or ``approximability") of tropical curves, whereas any harmonic tropical curve in $\mathbb{R}^m$ arise as the Hausdorff limit of families of harmonic amoebas. Obstructions to realisability in the algebraic setting are discussed for instance in \cite{Mikh06} or \cite{Nish}. There is at least one natural condition under which tropical curves are realisable. This condition, called regularity, is the one we will use in the remaining of the text. 

\subsection{Phase-tropical morphisms}\label{sec:ptm}
Recall that a toric morphism $ B : (\bC^\star)^k \rightarrow (\bC^\star)^m $ is a map of the form 
\[ (z_1, \dots, z_k) \mapsto \big( b_1 z_1^{a_{11}}\dots z_k^{a_{1k}}, \dots, b_m z_1^{a_{m1}}\dots z_k^{a_{mk}} \big) \]
where $(b_1,\dots,b_m) \in (\bC^\star)^m $ and $  \big( a_{ij} \big)_{ij} \in M_{m \times k} \big( \mathbb{Z} \big)$.

\begin{definition}\label{defphtropmorph}
A \df{phase-tropical morphism} $\phi : V \rightarrow (\bC^\star)^m $ on a simple phase-tropical curve $V = \big(\left\lbrace \pi_v, L_v \right\rbrace_{v\in V(C)}, \left\lbrace B_{\vec e} \right\rbrace_{e\in E(C)} \big)$ is a proper continuous map such that for any $v \in V(C)$, and the corresponding chart $\A^{-1}\big(\pi_v(T_v)\big)\cap L_v$, there exists a toric morphism $B_v : \ttor \rightarrow (\mathbb{C}^\star)^m$ such that the restriction of $\phi$ to the latter chart is given by $B_v$. We denote by $\pi_\phi : C \rightarrow \mathbb{R}^m$ the induced tropical morphism.
\end{definition}

\begin{remark}
Phase-tropical curves originally appeared as immersed objects in $\ttor$, obtained by degeneration of families of algebraic curves, see section $6.2$ in \cite{Mikh05}. Their abstract counterpart and associated morphisms were considered later in \cite{Mikhunpub}.
\end{remark}

\begin{proposition}\label{prop:phasetropdetermined}
Any phase-tropical morphism $\phi : V \rightarrow (\bC^\star)^m $ is uniquely determined by its underlying tropical morphism $ \pi_\phi : C \rightarrow \mathbb{R}^m$, up to a toric translation $(z_1, \dots, z_m) \mapsto \big( b_1 z_1, \dots, b_m z_m \big)$ on the target space.
\end{proposition}

\begin{proof}
Since the toric morphism $B_v : (z_1, \dots, z_k) \mapsto \big( b_1 z_1^{a_{11}}\dots z_k^{a_{1k}}, \dots, b_m z_1^{a_{m1}}\dots z_k^{a_{mk}} \big) $ descends to the affine-linear map $$ (x_1, \dots, x_k) \mapsto \big( \log\vert b_1\vert + a_{11} x_1+ \dots+ a_{1k} x_k,  \dots, \log \vert b_m \vert  + a_{m1} x_1 + \dots + a_{mk} x_k \big)$$ by composition with $\A$, the tropical morphism $\pi_\phi$ determines  $B_v$ up to the coordinate-wise argument of $b_v := (b_1,\dots,b_m) \in (\mathbb{C}^\star)^m $. Up to translation, we can assume that the latter vector is $(1,\dots,1)$ for a chosen vertex $v_0 \in V(C)$. For any edge $\vec e$ from $v_0$ to another vertex $v$, the argument of $b_v$ is determined by $B_{\vec e}$ and hence all the $b_v$ are determined up to a single translation by $\pi_\phi$. The result follows.
\end{proof}

Let us now recall the notions of regularity and superabundancy for tropical morphisms. We refer to \cite[\S 2.4-2.6]{Mikh05} for more details.
Recall that any simple tropical curve can be presented as a couple $(G,\ell)$ consisting of a cubic graph and a length function.

\begin{definition}
Let $C$ and $C'$ be two simple tropical curves supported on the same cubic graph $G$. Two tropical morphisms $\pi :  C \rightarrow \R^m$ and $\pi'  :  C' \rightarrow \R^m$ have the same \df{combinatorial type} if for any oriented leaf/edge $\vec{e} \in LE(G)$, 
$\pi(\vec{e})$ and $\pi'(\vec{e})$ have the same slope in $S^{m-1} \cup \left\lbrace 0 \right\rbrace$.
\end{definition}

The space of deformation of tropical morphisms within a fixed combinatorial type  is subject to a finite number of constraints. Therefore, its dimension is bounded from below by the ``expected dimension". Here, we are interested in the space of tropical curves supporting a tropical morphism of a fixed combinatorial type rather than the space of deformation of the morphism. We refer to 2.4 in \cite[\S 2.4]{Mikh05} for a proof of the next statement.

\begin{proposition}
Let $G$ be a cubic graph of genus $g$ with $n$ leaves and $R$ be a collection of residues of dimension $m$ on $G$. For any length function $\ell_0$ on $G$, the space of length function $\ell$ such that the tropical morphism $\pi_R : (G,\ell) \rightarrow \R^m$ has the same combinatorial type than $\pi_R : (G,\ell_0) \rightarrow \R^m$ is the relative interior of an open convex polyhedral domain of codimension at most $mg$ in the space $( \R_{>0})^{3g-3+n}$ of length functions on $G$.
\end{proposition}

The constraints on the above space of tropical curves are given by the $g$ cycles of $G$. Since the image of each of the $g$ cycles has to close up in $\R^m$, each cycle imposes $m$ conditions so that we obtain $mg$ conditions in total.

\begin{definition}
Let $G$ and $R$ be as above.  A tropical morphism $\pi_R : (G,\ell_0) \rightarrow \R^m$ is \df{regular} (respectively \df{superabundant}) if the space of length functions $\ell$ such that $\pi_R : (G,\ell) \rightarrow \R^m$ has the same combinatorial type than $\pi_R : (G,\ell_0) \rightarrow \R^m$ has codimension $mg$ (respectively smaller than $mg)$ in $(\R_{>0})^{3g-3+n} $. A phase-tropical morphism $\phi:V\rightarrow (\bC^\star)^m$ is \df{regular} if its underlying tropical morphism $\pi_\phi$ is.
\end{definition}
\medskip

\subsection{Mikhalkin's approximation Theorem}\label{Mik}

In this last section, we state and prove Mikhalkin's Theorem in a framework close to the one of \cite{Mikhunpub}. The main difference is the way regularity steps in the proof. In \cite{Mikhunpub}, regularity allows to find principal divisor on appropriate sequences of Riemann surfaces, giving the desired sequence of maps up to the Abel-Jacobi Theorem.

The approximation of phase-tropical morphisms requires more than the approximation of the underlying tropical morphism, guaranteed by Theorem \ref{approxtrop}. In the present framework, we need to construct sequences of harmonic amoeba maps that have well-defined harmonic conjugate and lift to actual holomorphic maps. Technically, we need harmonic amoeba maps coming from i.n.d. whose periods are integer multiples of $2\pi i$. The existence of such maps will be guaranteed by regularity, see Proposition \ref{propreg}.

Recall the change of complex structure 
\[
\begin{array}{rcl}
H_t \: : \: (\bC^\star)^m & \rightarrow & (\bC^\star)^m \\
(z_1,\dots,z_m) & \mapsto & \left( \vert z_1 \vert^{\frac{1}{\log(t)}} \frac{z_1}{\vert z_1 \vert} ,\dots, \vert z_m \vert^{\frac{1}{\log(t)}} \frac{z_m}{\vert z_m \vert} \right)
\end{array}.
\]
Then, we have the following formulation of Mikhalkin's approximation Theorem.

\begin{mytheorem}\label{thmMik}
Let $ \phi : V(G,\ell,\theta) \rightarrow (\bC^\star)^m$ be a regular phase-tropical morphism on a simple phase-tropical curve $V:=V(G,\ell,\theta)$, where $C:=(G,\ell)$ has genus $g$ and $n$ leaves. Then, there exists a sequence $\left\lbrace S_t \right\rbrace_{t>1} \subset \mgn$ converging to $V$, together with algebraic maps $ \phi_t \, : \, S_t \rightarrow (\bC^\star)^m$ such that $ H_t \big( \phi_t ( S_t ) \big) $ converges in Hausdorff distance to $\phi (V)$. Moreover, we can require the twist function $\theta_t$ of $S_t$ to be independent of $t$.
\end{mytheorem}

\begin{remark} 
--The requirement in Theorem \ref{thmMik} that $\theta_t$ does not depend on $t$ is of crucial importance for the construction of real algebraic curves with prescribed topology and in real enumerative geometry, see for instance \cite[Theorem $3$]{L} and \cite[Theorem $3$]{Mikh05}. Indeed, if $V$ is globally invariant under complex conjugation, then the twist function $\theta$ takes values in $\{-1, \,1\}$. If $\theta_t$ is chosen to be constant (and therefore equal to $\theta$), each Riemann surface $S_t$ inherits an anti-holomorphic involution given by the involution $\sigma$ of Theorem \ref{thm:buser} on every pair of pants of the underlying decomposition of $S_t$. Whenever the initial point of Definition \ref{defpermat}  in $S_t$ belongs to the fixed locus of the involution, the map $\phi_t$ is real algebraic.

-- Theorem \ref{thmMik} takes care of every phase-tropical immersion to $\ttor$ since any such morphism is regular, see \cite[Proposition $2.23$]{Mikh05}.
\end{remark}

\begin{proposition}\label{approxcomptrop}
Let $ \phi : V \rightarrow (\bC^\star)^m$ be a phase-tropical morphism on a simple phase-tropical curve $ V:=V(G,\ell,\theta)$ and $R$ be the collection of residues such that $\pi_\phi : C \rightarrow \mathbb{R}^m$ is equal to $\pi_R$. For any sequence $ \left\lbrace S_t \right\rbrace_{t>1} $ converging to $V$ and such that $\mathcal{P}_{R,S_t}$ is a constant family of integer period matrices, then $H_t \big(\phi_R(S_t)\big) $ converges in Hausdorff distance to $\phi(V)$, see Definition \ref{defpermat}.
\end{proposition}

\begin{proof}
We will proceed as in the proof of Theorem \ref{approxtrop}: we first show the convergence locally  and then check that the pieces glue together in the expected way.

Let us first show that $ H_t \big( \phi_{R} ( Y_{v,t} ) \big) $ converges to $\phi \big(  \A_V^{-1} (\htr _v) \big)$ for any $v \in V(C)$. Recall the notations $R_{v,t}$ and $R_v$ introduced in the proof of Theorem \ref{approxtrop}. By assumption, these matrices have integer coefficients and are all equal to each other. In particular, if $Y_{v,t}$ is presented as a subset of $S$ as in \eqref{eq:model}, the family of maps $\phi_R: Y_{v,t} \rightarrow (\bC^\star)^m$ converges point-wise to the map $\phi_{R_v}: S \rightarrow (\bC^\star)^m$ when a common initial point $p$ is chosen in $\cap_t Y_{v,t}\subset S$.
Recall also that $\phi_{R_v}: S \rightarrow (\bC^\star)^m$ is an embedding if and only if the $m\times 3$ matrix $R_v$ has rank $2$. If the latter map is not an embedding, consider the matrix $\widetilde R_v$ obtained by adding the matrix $R_0:=\begin{psmallmatrix} -1 & \;\;0 & \; 1\\ \;\;0 & -1 & \;1 \end{psmallmatrix}$ at the bottom of $R_v$ and consider the map $\widetilde \phi_{R} : Y_{v,t} \rightarrow (\bC^\star)^{m+2}$ obtained by adding the $2$ coordinates $e^{\int \omega_{R_0}^S}$ where the pair of i.n.d. $\omega_{R_0}^S$ restricts to $Y_{v,t}\subset S$. Since the map $\phi_{R}:Y_{v,t}\rightarrow (\bC^\star)^m$ (respectively $\phi_{R_v}: S \rightarrow (\bC^\star)^m$) is obtained from $\widetilde \phi_{R}$ (respectively $\phi_{\widetilde R_v}$) by projecting on the $m$ first coordinates, there is no loss of generality in assuming that $R_v$ has rank $2$, and then that  $\phi_{R_v}$ is an embedding. As observed above, the map $\phi_R: Y_{v,t} \rightarrow (\bC^\star)^m$ converges point-wise to $\phi_{R_v}: S \rightarrow (\bC^\star)^m$. We also know from Theorem \ref{approxtrop} that $\A\circ H_t\circ \phi_R(Y_{v,t})$ converges to $\pi_\phi(\htr_v)$. In particular, each boundary geodesic of $Y_{v,t}$ is mapped to the corresponding end of $\pi_\phi(\htr_v)$. From there, showing that $ H_t \big( \phi_{R} ( Y_{v,t} ) \big) $ converges to $\phi\big(\A_V^{-1} (\htr _v)\big)$ amounts to show that $H_t(\cL)$ converges to the phase-tropical line $L$. Observe in particular that the image of $\R Y_{v,t}$ converges to the real part of $\phi\big(\A_V^{-1} (\htr _v)\big)$.

It remains to show that for two vertices $v_1$, $v_2$ connected by an edge $e$,  the piece $H_t \big( \phi_{R} ( Y_{v_1,t} \cup Y_{v_2,t})\big)$ converges to $\phi\big(\A_V^{-1} (\htr _{v_1}\cup \htr _{v_2})\big)$. On the hyperbolic side, the two pairs of pants $Y_{v_1,t}$ and $Y_{v_2,t}$ are glued to each other according to the twist $\theta_t(e)$ measuring the angle between the two markings of $\gamma_{e,t}$ given by $\R Y_{v_1,t}$ and $\R Y_{v_2,t}$. On the phase-tropical side, the two charts $\A_V^{-1} (\htr_{v_1})$ and $\A_V^{-1} (\htr_{v_2})$ are glued to each other according to the twist $\theta(e)$ measuring the angle between the two markings of the geodesic $\A_V^{-1}(q)$ given by $\R \A_V^{-1} (\htr_{v_1})$ and $\R \A_V^{-1} (\htr_{v_1})$ where $q$ is the midpoint of $e$. The result follows now from the facts that: the twist $\theta_t(e)$ converges to $\theta(e)$, the image of $\R Y_{v_j,t}$ under $H_t\circ \phi_R$ converges to $\phi\big(\R \A_V^{-1} (\htr_{v_j})\big)$ and from Remark \ref{rem:blup}.
\end{proof}

\begin{lemma}\label{intperiod}
Let $ \phi : V \rightarrow (\bC^\star)^m$ be a phase-tropical morphism on a simple phase-tropical curve $ V:=V(C,\theta)$ and denote $R$ the collection of residues of the underlying tropical morphism $\pi_\phi : C \rightarrow \mathbb{R}^m$. For any loop $\rho\subset C$, we have 
\[  \sum_{\vec e \in \rho} \log \big( \theta (e) \big) \cdot \sw_{R}^{C} (\vec e) \in ( 2 i \pi \mathbb{Z} )^m. \]
\end{lemma}

\begin{proof}
For any vertex $v \in \rho$, there is a distinguished point $p_v$ among the 3 vertices of the coamoeba $\mathcal{A}_V^{-1} (v)$, namely the intersection point of the $2$ boundary geodesics of $\mathcal{A}_V^{-1} (v)$ corresponding to the 2 edges in $\rho$ adjacent to $v$. 
If $\vec e \subset \rho$ join the vertex $v$ to the vertex $v'$, the image of the distinguished point $p_{v'}$ under $\Arg \circ \phi $ is obtained from the image of $p_v$ by adding $\frac{1}{i}  \cdot \theta(e)  \cdot \sw_{R}^{C} (\vec e)$ in the argument torus $(\R/2\pi\bZ)^m$ of $(\bC^\star)^m$. We deduce that for any $v\in V(C)$, we have the equality 
\[ p_v + \sum_{\vec e \in \rho} \frac{1}{i} \log \big( \theta (e) \big) \cdot \sw_{R}^{C} (\vec e) \; = \; p_v\]
in the argument torus. The statement follows.
\end{proof}

\begin{lemma}\label{lem:convperiodinteger}
Let $ \phi : V \rightarrow (\bC^\star)^m$ be a phase-tropical morphism on a simple phase-tropical curve $ V:=V(C,\theta)$ and denote $R$ the collection of residues of the underlying tropical morphism $\pi_\phi : C \rightarrow \mathbb{R}^m$. For any family $ \left\lbrace S_t \right\rbrace_{t >1} \subset \mgn $ 
converging to $V$, the family of period matrices $ \left\lbrace \mathcal{P}_{R,S_t} \right\rbrace_{t >1} 
$ converges to an integer period matrix $ \mathcal{P}_{R, C}$, see Definition \ref{defpermat}.
\end{lemma}

\begin{proof}
According to Theorem \ref{convdif}, the period vector $ \frac{1}{2 i \pi} \big( \int_{\gamma_{e,t}} \omega_{R,1}^{S_t},\dots, \int_{\gamma_{e,t}} 
\omega_{R,m}^{S_t} \big)$ tends to the vector $\sw_{R}^{C} (\vec e)$, where $\gamma_{e,t}$ is 
the geodesic of the pair of pants decomposition of $S_t$ corresponding to $e \in E(C)$. Since $\pi_\phi$ is a tropical morphism, the vector $\sw_{R}^{C} (\vec e)$ has integer coordinates for any $e\in E(C)$. Hence, the limit of the above period vector has integer coordinates for any $e\in E(C)$.
Now, let us consider a basis $\rho_1,\dots, \rho_g \in H_1(C, \mathbb{Z})$ and consider the associated family of piecewise geodesic loops $\rho_{1,t},\dots, \rho_{g,t} \subset S_t$, see Definition \ref{assloop}. By Theorem \ref{thm:convperiod} and Lemma \ref{intperiod}, we have that 
\[\lim_{t\rightarrow \infty} \int_{\rho_{k,t}} \omega_{R,j}^{S_t} \in 2 i \pi \mathbb{Z}\]
for $1 \leqslant k \leqslant g$ and $1 \leqslant j \leqslant m$. Since the closed curves $\gamma_{t,e}$, $e\in E(C)$, and $\rho_{k,t}$, $1 \leqslant k \leqslant g$, generate $H_1(S_t,\mathbb{Z})$ for any $t$, the statement follows.
\end{proof}

For the rest of this section, the map $\phi : V \rightarrow (\bC^\star)^m$ is a regular phase-tropical morphism on a phase-tropical curve $ V:=V(C,\theta)$ and $R$ denote the collection of residues  of the underlying tropical morphism $\pi_\phi: C \rightarrow \mathbb{R}^m$. The genus and number of leaves of $C$ will be denoted $g$ and $n$ respectively and satisfy $2g+n>2$.
We denote by $G$ the cubic graph supporting $C$ and fix a ribbon structure $\scrR$ on it. Recall from Section \ref{sec:hyp} that we can define the surjective analytic map
\[
\begin{array}{rcl}
\fn_{G} \: : \: (\bC^\star)^{3g-3+n}  & \rightarrow & \mgn \\
(\ell,\theta) & \mapsto & S(G,\ell, \theta)
\end{array}.
\]
Consider the partial compactification $K$ of $(\bC^\star)^{3g-3+n}$ given by first applying the self- diffeomorphism 
\[(\ell,\theta)\mapsto \left(\frac{\ell}{\Vert \ell\Vert} (\Vert \ell\Vert+1), \theta \right)\]
and then taking the closure of the image in $\big(\R_{>0}\times (S^1)\big)^{3g-3+n} \subset (\bC^\star)^{3g-3+n}$. The latter construction corresponds to the real oriented blow-up of the $\ell$-factor of $(\bC^\star)^{3g-3+n}$ at the origin. The total space $K$ is diffeomorphic to $ \left\lbrace \ell \in (\R_{>0})^{3g-3+n} \: \big| \: \Vert \ell \Vert \geqslant 1 \right\rbrace\times (S^1)^{3g-3+n}$ with ``exceptional divisor"
$$E := \left\lbrace \ell \in (\R_{>0})^{3g-3+n} \: \big| \: \Vert \ell \Vert = 1 \right\rbrace \times   (S^1)^{3g-3+n}. $$ 
In view of the notion of abstract phase-tropical convergence given in Definition \ref{comptropconv}, it is natural to consider the points of $E$ as equivalence classes of phase-tropical curves with respect to the relation

\begin{equation}\label{eq:eqrelation}
V(G,\ell_1, \theta_1) \sim  V(G,\ell_2,\theta_2) \;\;  \Leftrightarrow \;\;  \big(\theta_1 = \theta_2 \; \text{ and } \; \ell_1 = \lambda \ell_2 \; \text{ for some } \; \lambda>0 \big).
\end{equation}
The map $\fn_{G} $ extends to $K$ and maps the exceptional divisor $E$ to the stable curve $S(G)$ dual to $G$. Now, define $k:=(2g+n-1)$ and consider the map 
\[
\begin{array}{rcl}
\Pi_R \: : \: (\bC^\star)^{3g-3+n}  & \rightarrow & M_{k \times m} (\mathbb{R}) \big/ \SL_{k} (\mathbb{Z}) \\
(\ell,\theta) & \mapsto & \mathcal{P}_{R,S(G,\ell,\theta)}
\end{array}
\]
where $\mathcal{P}_{R,S(G,\ell,\theta)}$ is the period matrix of $S(G,\ell,\theta)$ with respect to $R$, see Definition \ref{defpermat}.

\begin{proposition}\label{propreg}
The map $\Pi_R$ extends analytically to $K$. If $[V]$ denotes the equivalence class of $V$ in $E$, then each irreducible component of the level set $ \Pi_R^{-1} \big( \Pi_R( [V]) \big) \subset K$ is a smooth analytic subset of real codimension $2mg$ near $E$.
\end{proposition}

\begin{proof}
For any continuous family $\left\lbrace (\ell_t, \theta_t) \right\rbrace_{t>1}\subset (\bC^\star)^{3g-3+n}$ converging to a point $E\subset K$, the sequence $S_t:=\fn_G(\ell_t, \theta_t)$ converges to a phase-tropical curve $V$, up to an appropriate re-parametrisation of $t$. For any other re-parametrisation converging to a phase-tropical curve $V'$, the curves $V$ and $V'$ are equivalent under the relation \eqref{eq:eqrelation}. By Corollary \ref{cor:period}, the family of period matrices $\cP_{R,S_t}$ converges to a matrix that only depend on the equivalence class of $V$ and hence does not depend on the choice of a re-parametrisation. It follows that $\Pi_R$ extends to $K$. As seen in Corollary \ref{cor:period}, the limiting periods are polynomial in the coefficients of $\sw_{R}^{C}
(\vec e)$, $e\in E(C)$, and in the coefficients of $\theta$. Since the $\sw_{R}^{C}
(\vec e)$ depend linearly on $\ell$, the extension of $\Pi_R$ is analytic. 

According to Theorems \ref{thmcodim} and \ref{thmhyperb}, the level set $\Pi_R^{-1} \big( \Pi_R( [V]) \big)$ is the intersection of $m$ smooth subvarieties of real codimension $2g$ in $(\bC^\star)^{3g-3+n}$. In order to prove the statement, it remains to show that the level set of the map $\Pi_R$ extended to $K$ intersects $E$ in a smooth subvariety of codimension $2mg$. 
%
We claim that this intersection is described by the classes $[V']$ with $V'= V(C',\theta')$ and $C:=(G,\ell')$ satisfying
\begin{equation}\label{condperiod1}
\sum_{\vec e \in \rho} \ell'(e) \cdot  \sw_{R}^{C} (\vec e) = 0 \in \mathbb{R}^m
\end{equation}
\begin{equation}\label{condperiod2}
\sum_{\vec e \in \rho} \log \big( \theta'(e) \big)\cdot \sw_{R}^{C} (\vec e) = 0 \in (\R/2\pi\bZ)^m
\end{equation}
for any loop $\rho \subset G$. To see this, observe first that for any $\vec e\in E(G)$, we have $\sw_{R}^{C} (\vec e) =\sw_{R}^{C'} (\vec e)$. Indeed, by Theorem \ref{convdif}, the vector $\sw_{R}^{C} (\vec e)$ is the limit of $\int_{\gamma_{e,t}} \omega_{R}^{S_t}$ where $S_t:=S(G,\ell_t,\theta_t)$ and $(\ell_t,\theta_t)$ is any appropriately scaled family in $\Pi_R^{-1} \big( \Pi_R( [V]) \big)$ converging to $(\ell,\theta)$. The same holds for $\sw_{R}^{C'} (\vec e)$ and we have $\int_{\gamma_{e,t}} \omega_{R}^{S_t}=\int_{\gamma'_{e,t}} \omega_{R}^{S'_t}$ since both families $(\ell_t,\theta_t)$ and $(\ell'_t,\theta'_t)$ are in $\Pi_R^{-1} \big( \Pi_R( [V]) \big)$. It follows that $\sw_{R}^{C} (\vec e) =\sw_{R}^{C'} (\vec e)$. Now, the collection of equations \eqref{condperiod1} is equivalent to the exactness of $\sw_{R}^{C'} (\vec e)$ under the latter equalities. Therefore, these equations are satisfied on $\Pi_R^{-1} \big( \Pi_R( [V]) \big)\cap E$. For the collection of equations \eqref{condperiod2}, the left-hand side is, according to Theorem \ref{thm:convperiod} and the equalities $\sw_{R}^{C} (\vec e) =\sw_{R}^{C'} (\vec e)$, the limit of the period vector $\int_{\rho'_{t}} \omega_R^{S'_t}$ where $\rho'_t \subset S'_t$ is the loop associated to $\rho\subset C'$, see Definition \ref{assloop}. By the same argument as above, we have that $\int_{\rho'_{t}} \omega_R^{S'_t}=\int_{\rho_{t}} \omega_R^{S_t}$ and also that $\int_{\rho_{t}} \omega_R^{S_t}\in(2\pi i \bZ)^m$ by Theorem \ref{thm:convperiod} and Lemma \ref{intperiod}. It implies that the equations \eqref{condperiod2} are also satisfied on $\Pi_R^{-1} \big( \Pi_R( [V]) \big)\cap E$. In order to see that the equations \eqref{condperiod1} and \eqref{condperiod2}  are sufficient to describe $\Pi_R^{-1} \big( \Pi_R( [V]) \big)\cap E$, we need to show that the rank of the entire system of equations is $2mg$. Clearly, the system \eqref{condperiod1} is linearly independent from system \eqref{condperiod2} and the two systems have the same rank. Now, the system \eqref{condperiod1} describe the space of tropical curves $(G,\ell')$ for which the morphism $\pi_R$ has the same combinatorial type as $\pi_\phi$. Since $\pi_\phi$ is regular by assumption, the latter space as codimension $mg$. The statement follows.
\end{proof}

\begin{proof}[Proof of Theorem \ref{thmMik}.]
Recall that $2g+n>2$. Thanks to the above proposition, we can construct a sequence  $\left\lbrace (\ell_t,\theta_t) \right\rbrace_{ t>1} \subset \Pi_R^{-1} \big( \Pi_R ( [ V ] ) \big) $ converging to the point $\left[ V \right] \in E$. Up to a re-parametrisation, we can assume that $\ell_t \sim 2\pi^2 /\big( \log(t) \cdot \ell\big)$ where $\ell$ is the length function of $V$. In other words, the sequence $S_t : = \fn_{G} (\ell_t,\theta_t) $ converges to $V$ in the sense of  Definition \ref{comptropconv}. Since $\left\lbrace (\ell_t,\theta_t) \right\rbrace_{ t>1} \subset \Pi_R^{-1} \big( \Pi_R ( [ V ] ) \big) $,  the sequence of period matrices $ \cP_{R,S_t}$ is constant and integer, thanks to Lemma \ref{lem:convperiodinteger}. Applying Proposition \ref{approxcomptrop}, we deduces that $H_t \circ \phi_R (S_t) $ converges to $\phi(V)$.

To see that we can choose the sequence $ \left\lbrace \theta_t \right\rbrace_{t>1}$ to be constant, observe that the arguments of the above proof imply that $ \Pi_R^{-1} \big( \Pi_R ( [ V ] ) \big) \cap \left\lbrace \theta = \text{constant} \right\rbrace \subset K$ is analytic, smooth and of real codimension $mg+(3g-3+n)$ in a neighbourhood of $E$. In particular, the latter level set has strictly positive dimension. The statement follows.
\end{proof}

\bibliographystyle{alpha}
\bibliography{Draft}

\end{document}